\documentclass[a4,leqno,12pt]{article} 

\usepackage{amsmath, amsthm, amssymb} 

\numberwithin{equation}{section}

\title{ { \bf Higher Selberg Zeta Functions \\ for Congruence Subgroups } } 
\author{Tetsuya Momotani \\
\small
Graduate School of Mathematics, Kyushu University \\ 
\small 
6-10-1, Hakozaki Fukuoka 812-8581, Japan \\
\small 
E-mail: momo@math.kyushu-u.ac.jp
} 

\newtheorem{thm}{Theorem}[subsection]
\newtheorem{prop}[thm]{Proposition}
\newtheorem{lem}[thm]{Lemma}

\begin{document} 

\maketitle 

\begin{abstract}
As a generalization of the results \cite{Wakayama3}, 
we study the functional equation of the higher Selberg zeta function for congruence subgroups. 
To obtain the gamma factor of this function, 
we introduce a higher Dirichlet $L$-function. 
Then we determine the gamma factor explicitly 
in terms of the Barnes triple gamma function and the higher Dirichlet $L$-function. 
\end{abstract}

\section{Introduction}

For a discrete subgroup $ \varGamma $ of $ SL ( 2, \mathbb{R} ) $, 
the Ruelle zeta function is defined by 
$ \zeta _{ \varGamma }(s) := \prod_{ P } ( 1- N(P) ^{-s} )^{-1} $, 
where $ P $ runs through the all primitive hyperbolic conjugacy classes of $ \varGamma $, and $ N(P) $ is the norm of $P$. 
It is known (\cite{Ruelle}) 
that Selberg's zeta function $ Z_{ \varGamma } (s) $ 
is expressed as a shifted product of the Ruelle zeta function. 
In general, for a given zeta function $ z (s) $, 
the higher zeta function can be defined by the shifted product 
$ \prod z ( s+n ) $ (Cf. \cite{Wakayama1}). 
In this sense, $ Z_{\varGamma} (s) $ is the higher zeta function
constructed by the Ruelle zeta function. 

In \cite{Wakayama2}, the higher Selberg zeta function $ \prod _{ n=1 }^{ \infty } Z_{\varGamma } (s+n) ^{-1} $ was introduced 
for establishing a certain identity 
between the non-trivial zeros of Selberg zeta function and of the Riemann zeta function. 
This function also appears 
in the study \cite{Gon} of the first variations of the Selberg zeta function in Teichm\"{u}ller spaces. 
Further, analytic properties of the higher Selberg zeta function have already examined in \cite{Wakayama3} 
under the condition that the discrete subgroup $ \varGamma $ is a co-compact and torsion free. 
It was shown that this function also has a meromorphic continuation to the whole complex plane 
and satisfies a certain functional equation.

The purpose of this present paper is a generalization of the result in \cite{Wakayama3} to the non-compact cases. 
Precisely, we study the higher Selberg zeta function 
for the following congruence subgroup of $ SL ( 2, \mathbb{Z} ) $ 
for each integer $ N \geq 1 $: 
\begin{align*}
\varGamma_{0} (N) &:= \Big\{ 
\begin{pmatrix} a & b \\ c & d 
\end{pmatrix}
\in SL ( 2, \mathbb{Z} ) ; \; c \equiv 0 \, ( \bmod N ) \Big\} , \\
\varGamma_{1} (N) &:= \Big\{ 
\begin{pmatrix} a & b \\ c & d 
\end{pmatrix}
\in \varGamma_{0} (N) ; \; a \equiv d \equiv 1 \, ( \bmod N ) \Big\} , \\
\varGamma (N) \, &:= \Big\{ 
\begin{pmatrix} a & b \\ c & d 
\end{pmatrix}
\in \varGamma _{1} (N) ; \; b \equiv 0 \, ( \bmod N ) \Big\} .
\end{align*}

When $ \varGamma = \varGamma_{0} (N) $, $ \varGamma_{1} (N) $, $ \varGamma (N) $, 
the gamma factor of the Selberg zeta function $ \Xi _{ \mathrm{hyp} } ( s ) $ 
is composed of three factors $ \Xi _{ \mathrm{I} } ( s ) $, $ \Xi _{ \mathrm{ell} } ( s ) $, and $ \Xi _{ \mathrm{par} } ( s ) $, 
which respectively corresponds to central terms, elliptic terms, and parabolic terms of the trace formula. 
We show that these factors are expressed by the gamma function, 
the Barens double gamma function, and the Dirichlet $L$-function (Cf. \cite{Huxley},\cite{Vigneras}). 

\begin{align*}
\Xi _{ l \infty , \mathrm{hyp} } (s) 
:= \prod_{ m=1 }^{ \infty } \Xi _{ \mathrm{hyp} } ( s+lm )^{-1} .
\end{align*}
Noted that if $l=1$, 
this $ \Xi _{ l \infty , \mathrm{hyp} } (s) $ agrees with the one studied in \cite{Wakayama3} 
when $ \varGamma $ is a co-compact. 
In order to determine a proper gamma factor 
which describes the functional equation of $ \Xi _{ l \infty , \mathrm{hyp} } (s) $, 
we introduce a higher Dirichlet $L$-function. 
We explicitly determine three factors 
$ \Xi _{ l \infty , \mathrm{I} } (s) $, $ \Xi _{ l \infty , \mathrm{ell} } (s) $, $ \Xi _{ l \infty , \mathrm{par} } (s) $. 
We also describe the functional equation of the higher Dirichlet $L$-function. 

\section{Preliminaries} 

To investigate the higher Selberg zeta function, 
we briefly review the trace formula and the Selberg zeta function for congruence subgroups. 

\subsection{Selberg's Trace Formula} 

We first recall the trace formula. 
Let $ \varGamma = \varGamma_{0} (N), \varGamma_{1} (N), \varGamma (N) $ be the congruence subgroup 
and $ H = \{ z = x + iy ; \; x \in \mathbb{R} , y > 0 \} $ be the upper half plane. 
The group $ \varGamma $ acts discontinuously on $ H $ by linear fractional transformations. 
Let $ \Delta := - y^{2} \big( \frac{ d^{2} }{ dx^{2} } + \frac{ d^{2} }{ dy^{2} } \big)$ be the Laplacian on $ H $, 
which has a unique extension as a self adjoint operator 
acting on the space $ L^{2} ( \varGamma \setminus H ) $. 
We put the eigenvalues of 
$ \Delta $ on $ L^{2} ( \varGamma \setminus H ) $ 
by $ 0 = \lambda_{0} < \lambda_{1} \leq \lambda_{2} \leq \lambda_{3} \leq \cdots $, 
and $ r_{n} := \sqrt{ \lambda_{n} - 1/4 } $. 

When $ \varGamma = \varGamma_{0} (N), \varGamma_{1} (N), \varGamma (N) $, 
the trace formula contains the elliptic terms and the parabolic terms. 
To write down the parabolic terms, we need the description of the scattering matrix $ \varphi (s) $. 
It is known that (Cf. \cite{Selberg}) the determinant of the scattering matrix for $ \varGamma = SL ( 2, \mathbb{Z} ) $ is given by
\begin{align*}
\phi (s) 
= \det \varphi (s) 
= \frac{ \pi^{ -(1-s) } \Gamma ( 1-s ) \zeta ( 2-2s ) }{ \pi ^{-s} \Gamma (s) \zeta (2s) } ,
\end{align*}
where $ \zeta (s) $ is the Riemann zeta function. 
For $ \varGamma = \varGamma_{0} (N), \varGamma_{1} (N), \varGamma (N) $, 
Huxley (\cite{Huxley}) calculated the determinant of the scattering matrix as follows. 

\begin{lem}[\cite{Huxley}]
The determinant of the scattering matrix is given by
\begin{align}
\phi (s) = (-1)^{ ( \kappa -  \kappa _{0} ) / 2 } \Big( \frac{ \Gamma (1-s) }{ \Gamma (s)} \Big) ^{ \kappa } 
\Big( \frac{ \mathcal{A} }{ \pi ^{ \kappa } } \Big) ^{ 1-2s }
\prod_{ \chi } \frac{ L( 2-2s, \bar{ \chi } ) }{ L(2s, \chi) } . \label{ScaMat}
\end{align}
Here $ \kappa $ is the number of cusps and $\kappa _{ 0 } := - \mathrm{tr} \; \varphi \big( \frac{1}{2} \big) $. 
Dirichlet characters $ \chi $ which appear in the product of $(\ref{ScaMat})$ are expressed as 
$ \chi (n) = \chi_{1}(n) \, \chi_{2} (n) \, \omega _{ m_{1} m_{2} } (n) $,
where $ \chi_{i} $ is the primitive Dirichlet character modulo $ q_{i} $ $ ( i = 1, 2 ) $, and
$ \omega _{ m_{1} m_{2} }$ is the principal character modulo $ m_{1} m_{2} $. 
For each congruence subgroup $ \varGamma = \varGamma_{0} (N), \varGamma_{1} (N), \varGamma (N) $, 
$ \chi $ runs through all the pairs 
$ ( \chi_{1} , \chi_{2}, q_{1}, q_{2}, m_{1} ,m_{2} ) $ which satisfy the following conditions: 
\begin{align*}
\varGamma _{0} (N) \; &; \quad
\chi_{1} = \chi_{2} , \; q_{1} = q_{2} , \; m_{1} = 1 , \; q_{1} | m_{2} , \; ( m_{2}q_{2} ) | N , \; ( m_{1} , m_{2} ) = 1 , \\
\varGamma _{1} (N) \; &; \quad m_{1} = 1 , \; q_{1} | m_{2} , \; ( m_{2}q_{2} ) | N , \; ( m_{1} , m_{2} ) = 1 , \\
\varGamma (N) \, \; \, &; \quad ( m_{1}q_{1} ) | N, \; ( m_{2}q_{2} ) | N , \; ( m_{1} , m_{2} ) = 1 .
\end{align*}
A positive constant $ \mathcal{A} $ is expressed as 
\begin{align*}
\mathcal{A} = 
\begin{cases} 
\displaystyle \prod_{\,} \frac{ q_{1} N }{ (m_{1} , N/ m_{1} ) } & \text{on } \varGamma_{0} (N) , \\
\displaystyle \prod_{\,} q_{1} N & \text{on } \varGamma_{1} (N) , \\
\displaystyle \prod m_{1} m_{2} \,q_{1} N  & \text{on } \varGamma (N) , 
\end{cases}
\end{align*}
where $ ( \chi_{1}, \chi_{2},  q_{1}, q_{2}, m_{1} ,m_{2} ) $ 
runs over through all the pairs fulfilling above conditions. \qed
\end{lem} 

By making use of this lemma, we can write down the trace formula for the congruence subgroup as follows. 

\begin{thm}[\cite{Huxley}]
Suppose that the function $ h(r) $ is even, 
holomorphic, and $ h(r) = O ( (1 + |r|)^{ - 2 - \delta } ) $ 
in the strip $ | \mathrm{Im} (r) | \leq 1/2 + \delta $ $ ( \exists \delta > 0 ) $. 
Then we have 
\begin{align*}
\sum_{ n = 0 }^{ \infty } h (r_n) 
&= \frac{ \mathrm{vol} ( \varGamma \setminus H ) }{ 4 \pi } \int_{ - \infty } ^{ \infty } h(r) r \tanh ( \pi r) dr \\
&+ \sum_{ \{ P \} _{ \varGamma } } \sum_{k=1}^{ \infty } 
\frac{ \log N(P) }{ N(P)^{ \frac{k}{2} } - N(P)^{ - \frac{k}{2} } } g( k \log N(P) ) \\
&+ \frac{ \nu_{2} }{4} \int_{ - \infty }^{ \infty } h(r) \frac{1}{ e^{ \pi r } + e^{ - \pi r } }  dr 
+ \frac{ \nu_{3}  }{ 3 \sqrt{3} } \int_{ - \infty }^{ \infty } h(r) 
\frac{ e^{ \frac{ \pi r }{3} } + e^{ - \frac{ \pi r }{3} } }{ e^{ \pi r} + e^{ - \pi r } } dr \\
&- g(0) \log \big( \mathcal{A} \frac{ 2^{ \kappa } }{ \pi ^{ \kappa } } \big) 
+ \frac{1}{4} ( \kappa - \kappa _{0} ) h(0) \\
&- \frac{ \kappa }{ 2 \pi } \int_{ - \infty }^{ \infty } h(r) 
\big\{ \frac{ \Gamma ^{ \prime } }{ \Gamma }  (1 + ir) +  \frac{ \Gamma ^{ \prime } }{ \Gamma } ( \frac{1}{2}  + ir ) \big\} dr 
+ 2 \sum_{ \chi } \sum_{ n=1 }^{ \infty } \frac{ \chi (n) \Lambda (n) }{n} g(2 \log n) .
\end{align*}
Here, $ \mathrm{vol} ( \varGamma \setminus H ) $ is the volume of the fundamental domain of $ \varGamma $.  
$ \{ P \} _{\varGamma} $ runs through all the primitive hyperbolic conjugacy classes of $ \varGamma $, 
and $ N(P) := \mathrm{ max } \{ \alpha_{P} ^{2} , \beta_{P} ^{2} \} $, 
where $ \alpha _{P} $ and $ \beta _{P} $ are eigenvalues of the matrix $P$. 
The function $ g(u) $ is the inverse of the Fourier transform of $ h(r) $: 
\begin{align*}
g(u) = \frac{1}{ 2 \pi } \int_{ - \infty }^{ \infty } h(r) e^{ -iur } dr.
\end{align*}
Moreover, $ \nu_{2} $ and $ \nu_{3} $ is the number of primitive elliptic classes with order 2 and 3 respectively, 
and $ \Lambda (n)$ denotes the von Mangoldt function. \qed
\end{thm}

\subsection{Selberg Zeta Function}

We next recall several properties of Selberg zeta function. 
The Selberg zeta function for $ \varGamma $ is defined by the Euler product 
\begin{align*}
\Xi_{ \mathrm{hyp} }(s) 
:= \prod_{ n=0 }^{ \infty } \prod_{ \{ P \} _{ \varGamma } } \big( 1 -  N(P)^{-s-n} \big) , 
\end{align*}
where $ \{ P \} _{ \varGamma } $ runs through all the primitive hyperbolic conjugacy classes of $ \varGamma $. 
This product converges absolutely in $ \mathrm{Re} (s) > 1 $. 
We note that the logarithmic derivative of $ \Xi_{ \mathrm{hyp} }(s) $ is given by
\begin{align}
\frac{d}{ds} \log \Xi_{ \mathrm{hyp} }(s) = \sum_{ \{ P \} _{ \varGamma } } \sum_{ k=1 }^{ \infty } 
\frac{ \log N(P) }{ N(P)^{ \frac{k}{2} } - N(P)^{ - \frac{k}{2} } } \cdot N (P) ^{ - ( s - \frac{1}{2} ) k } .
\label{LDSZ}
\end{align}
By taking the test function
\begin{align*}
g(u) = \frac{1}{2s-1} e^{-(s-\frac{1}{2})|u|} - \frac{1}{2a-1} e^{-(a-\frac{1}{2})|u|} , 
\quad \mathrm{Re} (s)  > 1, a > 1 ,
\end{align*}
we obtain the analytic continuation and the functional equation of $ \Xi_{ \mathrm{hyp} }(s) $ as follows. 

\begin{thm}[Cf. \cite{Fischer}, \cite{Koyama}]
The Selberg zeta function $ \Xi _{ \mathrm{hyp} } (s) $ defined for $ \mathrm{Re} (s) > 1 $ 
has a meromorphic continuation to the whole complex plane, 
and the complete Selberg zeta function 
$ \Xi (s) : = \Xi_{\mathrm{I} }(s) \cdot \Xi_{\mathrm{hyp} }(s) \cdot \Xi_{\mathrm{ell} }(s) \cdot \Xi_{\mathrm{par} }(s) $
satisfies the functional equation 
\begin{align*}
\Xi (s) = \Xi (1-s) . 
\end{align*}
Here three factors $ \Xi_{\mathrm{I} }(s) $, $ \Xi_{\mathrm{ell} }(s) $, $ \Xi_{\mathrm{par} }(s)$ are explicitly given by 
\begin{align}
& \Xi_{ \mathrm{I} }(s) := \exp \Big\{ \frac{ \mathrm{vol} ( \varGamma \setminus H ) }{ 2 \pi }
\Big( s \log 2 \pi + \log \frac{ \Gamma _{2} (s) ^{2} }{ \Gamma (s) } \Big) \Big\} 
, \label{SelZetaGammaFactor} \\
& \Xi_{ \mathrm{ell} }(s) := 
\Big\{ \Gamma ( \frac{s}{2} ) ^{ -1 } \Gamma( \frac{s+1}{2} ) \Big\} ^{ \nu_{2} /2 }
\cdot \Big\{ \Gamma ( \frac{s}{3} ) ^{ -1 } \Gamma( \frac{s+2}{3} ) \Big\} ^{ 2 \nu_{3} /3 } , \notag \\
& \Xi_{ \mathrm{par} }(s) := 
\Big( \mathcal{A} \frac{ 2^{ \kappa } }{ \pi ^{\kappa} } \Big) ^{ - s} 
\Big( s - \frac{1}{2} \Big) ^{ ( \kappa - \kappa _{0} ) / 2 } \,
\Gamma ( s + \frac{1}{2} ) ^{ - \kappa } \, \Gamma ( s ) ^{ - \kappa } \, 
\prod_{ \chi } L ( 2s , \chi ) ^{-1} , \notag
\end{align}
and $ \Gamma _{2} (z) $ is the Barnes-Vign\'{e}ras double gamma function (Cf. \cite{Vigneras}) given by 
\begin{align}
\frac{1}{ \Gamma_{2} (z+1) } 
= \exp \Big\{ \big( - \frac{1}{2} -\frac{ \gamma }{2}  \big) z^{2} + \big( - \frac{1}{2} - \zeta ^{\prime} (0)  \big) z  \Big\} 
\prod _{k=1}^{\infty} \big( 1 + \frac{z}{k} \big) ^{k} \exp \big( -z + \frac{ z^{2} }{ 2k } \big) ,
\end{align}
which satisfies $ \Gamma_{2} (z+1) = \Gamma (z) ^{-1} \cdot \Gamma_{2} (z) $, and $ \Gamma_{2} (1) = 1 $.
Furthermore, the function $ \Xi (s) $ is an entire function of order $2$ with zeros at $ s = 1/2 \pm i r_{n} $ $ ( n \geq 0 ) $ only. 
\qed
\end{thm}

When $ \varGamma = SL( 2 , \mathbb{Z} ) $ is the modular group, 
it is known that
\begin{align*}
\mathrm{vol} ( \varGamma \setminus H ) = \frac{ \pi }{3} , \quad \nu_{2} = \nu_{3} =1 ,\quad  \mathcal{A} = 1, \quad \kappa = \kappa_{0} =1, 
\end{align*}
and $ L (s, \chi_{0} ) = \zeta (s) $ is the Riemann zeta function. 
Thus, poles and zeros of the Selberg zeta function 
$ \Xi_{\mathrm{hyp} }(s) = \Xi (s) \cdot \Xi_{\mathrm{I} }(s)^{-1} \cdot \Xi_{\mathrm{ell} }(s)^{-1} \cdot \Xi_{\mathrm{par} }(s)^{-1} $ 
for $ \varGamma = SL( 2, \mathbb{Z}) $ are explicitly given as follows. \\

\noindent{\bf Poles of $ \Xi_{\mathrm{hyp} }(s) $.}

(1) $ s=0 $ ; \; order $1$,

(2) $ s= 1/2 - k $ $( k \geq 0)$ ; \; order $1$. \\

\noindent{\bf Zeros of $ \Xi_{\mathrm{hyp} }(s) $.}

(1) $ s=1 $ ; \quad order $1$,

(2) $ s = - 6k-j $ $ (k \geq 0 \; , j=1,2,3,4,6 ) $ ; \; order $ 2k+1 $,

$ \;\; \quad  s= - 6k- 5 $ $ (k \geq 0 ) $ ; \; order $ 2k+3 $,

(3) $ s = \rho / 2 $ ( $ \rho $ : non-trivial zeros of $ \zeta (s) $) ,

(4) $ s = 1/2 \pm i r_{n} $ $ ( n \geq 1 ) $.

\section{Higher Selberg Zeta Functions}

In this section, 
we determine the gamma factor of the higher Selberg zeta function and describe the functional equation.

\subsection{Euler Product and Analytic Continuation}

We first define the higher Selberg zeta function attached to a positive integer $l \geq 1$. 
Let $ \Xi _{ \mathrm{hyp} } (s) $ be the Selberg zeta function. 
For each $ l \geq 1 $, 
we define the higher Selberg zeta function $ \Xi _{ l \infty , \mathrm{hyp} }(s) $ by the product 
\begin{align}
\Xi _{ l \infty , \mathrm{hyp} }(s) := \prod_{ m = 1 }^{ \infty } \Xi _{ \mathrm{ hyp } } ( s+lm ) ^{-1} 
= \prod_{m=1}^{\infty} \prod_{ n=0 }^{ \infty } \prod_{ \{ P \} _{ \varGamma } } \big( 1 - N(P)^{ -s - lm - n } \big) ^{-1}  .
\end{align}
This product converges absolutely in $ \mathrm{Re} (s) > 1-l $. 
Furthermore, the relation 
\begin{align}
\Xi _{ l \infty , \mathrm{hyp} } (s) 
= \frac{ \Xi _{ l \infty , \mathrm{hyp} } ( s+l ) }
{ \Xi _{ \mathrm{hyp} } (s+l) }  \label{HSZO}
\end{align}
shows that $ \Xi _{ l \infty , \mathrm{hyp} }(s) $ has a meromorphic continuation to the whole complex plane.

We note that the Euler product of the higher Selberg zeta function is rewritten as 
\begin{align*}
\Xi _{ l \infty , \mathrm{hyp} }(s) & = \prod_{ k = 1 }^{ \infty } \Xi _{ \mathrm{ hyp } } ( s+lk ) ^{-1} 
= \prod_{ k=0 }^{ \infty } \prod_{ n=0 }^{ \infty } \prod_{ \{ P \} _{ \varGamma } } 
\big( 1 - N(P) ^{ -s-l-lk-n } \big) ^{-1} \\
& = \prod_{ r=0 }^{ l-1 } \prod_{ k=0 }^{ \infty } \prod_{ j =0 }^{ \infty } \prod_{ \{ P \} _{ \varGamma } } 
\big( 1 - N(P) ^{ -s-l-lk-lj-r } \big) ^{-1} \quad ( n = lj+r ) \\
& = \prod_{ r=0 }^{ l-1 } \prod_{ m=1 }^{ \infty } \prod_{ \{ P \} _{ \varGamma } } 
\big( 1 - N(P) ^{ -s-lm-r } \big) ^{-m} \quad ( 1+k+j = m ). 
\end{align*}

\subsection{Test Function}

To determine the gamma factor of $ \Xi _{ l \infty , \mathrm{hyp} } (s) $, 
we need an appropriate test function. 
The following test function is essentilally the same one of \cite{Wakayama3}. 

\begin{prop}
If 
\begin{align}
g (u) = \frac{ u e^{ - ( s + \frac{l}{2} - \frac{1}{2} ) |u| } }{ 2 \sinh \frac{ lu }{2} } 
- \frac{1}{l} e^{ - ( s + \frac{l}{2} - \frac{1}{2} ) |u| }, \quad \mathrm{Re} (s) > 1-l , \label{TestFun}
\end{align}
then the Fourier transform of $ g(u) $ is given as 
$ h(r) = 
\phi _{s} (r)  + \phi _{s} ( -r ) $, 
where
\begin{align*}
\phi _{s} (r) 
:= & \sum_{ m=0 }^{ \infty } \Big\{ \frac{1}{ (s+lm+l-\frac{1}{2} +ir )^{2} } 
- \frac{1}{l} \frac{1}{ s+lm + \frac{l}{2} - \frac{1}{2} +ir }
+ \frac{1}{l} \frac{1}{ s+lm + \frac{ 3l }{2} - \frac{1}{2} +ir } \Big\} \notag \\
= &  - \frac{ l^{2} }{4} \sum_{ m=0 }^{ \infty } \frac{1}{ (s+lm+l-\frac{1}{2} +ir )^{2}  
( s+lm + \frac{l}{2} - \frac{1}{2} +ir ) ( s+lm + \frac{ 3l }{2} - \frac{1}{2} +ir )} .
\end{align*}
Moreover, this function satisfies the condition of the trace formula. 
\end{prop} 

\begin{proof} 
Similar to \cite{Wakayama3}. 
\end{proof}

Using this test function. we describe the hyperbolic terms of the trace formula as follows. 

\begin{prop}
Applying the function in $(\ref{TestFun})$, we have
\begin{align}
\sum_{ \{ P \} _{ \varGamma } } \sum_{ k=1}^{ \infty } 
\frac{ \log N(P) }{ N(P)^{ \frac{k}{2} } - N(P)^{ - \frac{k}{2} } } g ( k \log N(P) )
= \frac{ d^{2} }{ ds^{2} } \log \Xi _{ l \infty , \mathrm{hyp} } (s) 
- \frac{1}{l} \frac{d}{ ds } \log \Xi _{ \mathrm{hyp} } ( s + \frac{l}{2} ) ,\label{HC}
\end{align}
for $ Re(s) > 1-l $ . 
\end{prop}

\begin{proof} 
This equation follows immediately from the relation (\ref{LDSZ}) and 
\begin{align*}
\frac{ d^{2} }{ ds^{2} } \log \Xi _{ l \infty , \mathrm{hyp} } (s)  
= \sum_{ \{ P \} _{ \varGamma } } \sum_{ k=1}^{ \infty } \frac{ \log N(P) }{ N(P)^{ \frac{k}{2} } - N(P)^{ - \frac{k}{2} } } 
\frac{ k \log N(P) }{ N(P) ^{ \frac{ lk }{2} }  - N(P) ^{ - \frac{lk}{2} } } N (P) ^{ - ( s + \frac{l}{2} - \frac{1}{2} ) k } .
\end{align*}

\end{proof}

\subsection{Central Factor}

We calculate the central terms of the trace formula by taking the above mentioned test function. 

\begin{lem}
If $ h(r) = \phi _{s} (r)  + \phi _{s} ( -r ) $, we have
\begin{align}
\frac{1}{2} \int_{ - \infty }^{ \infty } h(r) r \tanh ( \pi r ) dr
& = \frac{ d^{2} }{ ds^{2} } \log \prod_{ r=0 }^{ l-1 } 
\Gamma _{3} \big( \frac{s+l+r}{l} \big) ^{ -2l } \Gamma _{2} \big( \frac{s+l+r}{l} \big) ^{ -2r + 2l -1 } \label{Lem311G} \\
& - \frac{1}{l} \frac{d}{ ds } \log \prod_{ r=0 }^{ l-1 }
\Gamma _{2} \big( \frac{s+ \frac{l}{2} +r}{l} \big) ^{ 2l } \Gamma \big( \frac{s+ \frac{l}{2} +r}{l} \big) ^{ 2r - 2l +1 }. \notag
\end{align}
Here $ \Gamma_{3} (z) $ is the Barnes-Vign\'{e}ras triple gamma function (Cf. \cite{Vigneras}) given by 
\begin{align} 
\Gamma_{3} (z+1) 
=& \exp \Big\{ 
\big( - \frac{1}{4} -\frac{ \gamma }{6} \big) z^{3} + \big( \frac{1}{8} - \frac{ \zeta ^{\prime} (0) }{2} + \frac{ \gamma }{4} \big) z^{2}  
+ \big( \frac{7}{24} + \zeta ^{\prime} (-1) + \frac{ \zeta ^{\prime} (0) }{2} \big) z  \Big\} \\ 
\times & \prod _{k=1}^{\infty} \Big( 1 + \frac{z}{k} \Big) ^{ - \frac{k(k+1)}{2} } 
\exp \Big\{ \frac{1}{6k} z^{3} + \big( - \frac{1}{4} - \frac{1}{4k} \big) z^{2} + \big( \frac{1}{2} +\frac{k}{2} \big) z \Big\} , \notag
\end{align} 
which satisfies $ \Gamma_{3} (z+1) = \Gamma_{2} (z) ^{-1} \cdot \Gamma_{3} (z) $, and $ \Gamma_{3} (1) = 1 $.
\end{lem}

\begin{proof}
Applying the residue theorem to the lower half plane, we see that
\begin{align*}
& \frac{1}{2} \int_{ - \infty }^{ \infty } h(r) r \tanh ( \pi r ) dr = \int_{ - \infty }^{ \infty } \phi_{s} (r) r \tanh ( \pi r ) dr 
= - \sum_{ n=0 }^{ \infty } ( 2n+1 ) \phi_{s} \big( - \frac{ 2n+1 }{2} i \big) \\
= & - \sum_{ n=0 }^{ \infty } \sum_{ m=0 }^{ \infty } ( 2n+1 ) 
\Big\{ \frac{1}{ (s+lm+l+n )^{2} } - \frac{1}{l} \frac{1}{ s+lm + \frac{l}{2} +n } 
+ \frac{1}{l} \frac{1}{ s+lm + \frac{3l}{2} +n } \Big\} \\
= & - \sum_{ r=0 }^{ l-1 } \sum_{ j=0 }^{ \infty } \sum_{ m=0 }^{ \infty } \big\{ 2(lj+r) +1 \big\} \\
& \qquad \times \Big\{ \frac{1}{ (s+lm+l+lj+r )^{2} } - \frac{1}{l} \frac{1}{ s+lm+ \frac{l}{2} +lj+r } 
+ \frac{1}{l} \frac{1}{ s+lm + \frac{ 3l }{2} +lj +r } \Big\} \\
= & - \sum_{ r=0 }^{ l-1 } \sum_{ k=0 }^{ \infty } \sum_{ j=0 }^{k} \big\{ 2(lj+r) +1 \big\} 
\Big\{ \frac{1}{ (s+lk+l+r )^{2} } - \frac{1}{l} \frac{1}{ s+lk+ \frac{l}{2} +r } 
+ \frac{1}{l} \frac{1}{ s+lk + \frac{ 3l }{2} +r } \Big\} \\
= & - \sum_{ r=0 }^{ l-1 } \sum_{ k=0 }^{ \infty } \Big\{ 2l \frac{ (k+1)(k+2) }{2} + (2r-2l+1)(k+1) \Big\} \\
& \qquad \times \Big\{ \frac{1}{ (s+lk+l+r )^{2} } - \frac{1}{l} \frac{1}{ s+lk+ \frac{l}{2} +r } 
+ \frac{1}{l} \frac{1}{ s+lk + \frac{ 3l }{2} +r } \Big\} .
\end{align*}
Here, from the Weierstrass product expression of $ \Gamma (z) $, $ \Gamma _{2} (z) $ and $ \Gamma_{3} (z) $, we have
\begin{align*}
& \frac{ d^{2} }{ ds^{2} } \log \Gamma _{3} \big( \frac{s+l+r}{l} \big) ^{2l} \Gamma _{2} \big( \frac{s+l+r}{l} \big) ^{2r-2l+1} \\
= & \frac{1}{ l^{2} } ( -3-2 \gamma ) s + \frac{1}{ l^{2} } ( \gamma -r+1 ) 
+ \frac{1}{l} ( - \gamma -2 \zeta ^{ \prime } (0) - \frac{3}{2} ) \notag \\
& +  \sum_{ k=0 }^{ \infty } \Big\{ \frac{2l \frac{ (k+1)(k+2) }{2} + (2r-2l+1)(k+1) }{ (s+lk+l+r )^{2} } 
+ \frac{1}{ l^{2} } \frac{ 2s+l-1 }{ k+1 } - \frac{1}{l} \Big\} , 
\end{align*}
and
\begin{align*}
& \frac{d}{ ds } \log 
\Gamma _{2} \big( \frac{ s+ \frac{l}{2} +r }{l} \big) ^{ 2l } \Gamma \big( \frac{ s + \frac{l}{2} + r }{l} \big) ^{ 2r-2l+1 } \\
= & \frac{1}{l} ( 2 + 2 \gamma ) s + \frac{1}{l} ( 2r - \gamma ) + \gamma + 2 \zeta ^{ \prime } (0) \\
& -  \sum_{ k=0 }^{ \infty } \Big\{ \frac{ 2l (k+1) + 2r-2l+1 }{ s+lk+ \frac{l}{2} +r } 
+ \frac{ 2s+l-1 }{ l (k+1) } - 2 \Big\} . \qquad \qquad \quad \quad
\end{align*}
It follow that
\begin{align*}
& \sum_{ r=0 }^{ l-1 } \sum_{ k=0 }^{ \infty } \Big\{ 2l \frac{ (k+1)(k+2) }{2} + (2r-2l+1)(k+1) \Big\} \\
& \qquad \times \Big\{ \frac{1}{ (s+lk+l+r )^{2} } - \frac{1}{l} \frac{1}{ s+lk+ \frac{l}{2} +r } 
+ \frac{1}{l} \frac{1}{ s+lk + \frac{ 3l }{2} +r } \Big\} \\
- & \frac{ d^{2} }{ ds^{2} } 
\log \Gamma _{3} \big( \frac{s+l+r}{l} \big) ^{2l} \Gamma _{2} \big( \frac{s+l+r}{l} \big) ^{2r-2l+1} 
- \frac{1}{l} \frac{d}{ ds } \log 
\Gamma _{2} \big( \frac{ s+ \frac{l}{2} +r }{l} \big) ^{ 2l } \Gamma \big( \frac{ s + \frac{l}{2} + r }{l} \big) ^{ 2r-2l+1 } \\
= & \frac{s}{ l^{2} } - \frac{1}{ l^{2} } ( r+1 ) + \frac{3}{2l} -  \frac{1}{l} \sum_{ k=0 }^{ \infty } 
\frac{ ( s+ \frac{l}{2}+r ) ( s+ \frac{3l}{2} -r-1 ) } { (s+lk+ \frac{l}{2} +r )(s+lk+ \frac{3l}{2} +r ) } \\
= & \frac{s}{ l^{2} } - \frac{1}{ l^{2} } ( r+1 ) + \frac{3}{2l} 
-  \frac{1}{ l^{2} } ( s+ \frac{l}{2}+r ) ( s+ \frac{3l}{2} -r-1 ) 
\sum_{ k=0 }^{ \infty } \Big\{ \frac{1} { s+lk+ \frac{l}{2} +r } - \frac{1}{ s+lk+ \frac{3l}{2} +r } \Big\} \\
= & \frac{s}{ l^{2} } - \frac{1}{ l^{2} } ( r+1 ) + \frac{3}{2l} - \frac{1}{ l^{2} } ( s+ \frac{3l}{2} -r-1 ) = 0 .
\end{align*}
This proves the lemma. 
\end{proof}

Now we recall the multiplication formula for $ \Gamma (z) $ and $ \Gamma _{2} (z) $ (Cf. \cite{Srivastava}). 
\begin{align}
& \Gamma (z) = k_{1} (m) \cdot m^{z} \cdot \prod_{ j=0 }^{ m-1 } \Gamma ( \frac{ z+j }{m} ) , \label{GaMult} \\
& \Gamma_{2} (z) = k_{2}(m) \cdot ( 2 \pi )^{ \frac{1}{2} (m-1) z } \cdot m^{ z- \frac{1}{2} z^{2} } 
\cdot \prod_{ j_{1} = 0 }^{ m-1 } \prod_{ j_{2} = 0 }^{ m-1 } \Gamma_{2} \big( \frac{ z + j_{1} + j_{2} }{m} \big) , \label{DGMult}
\end{align} 
where 
\begin{align*} 
& k_{1} (m) 
:= ( 2 \pi )^{ - \frac{1}{2} (m-1) } \cdot m^{- \frac{1}{2} } ,\\
& k_{2} (m) 
:= A ^{ 1 - m^{2} } \cdot e^{ \frac{1}{12} ( m^{2} -1 ) } 
\cdot ( 2 \pi )^{ - \frac{1}{2} (m-1) } \cdot m^{ - \frac{5}{12} }, 
\end{align*} 
and $A$ is the Glaisher-Kinkelin constant defined by
\begin{align*} 
\log A := \lim _{ N \to \infty } \Big\{
\sum_{k=1}^{N} k \log k - \Big( \frac{ N^{2} }{2} + \frac{N}{2} + \frac{1}{12} \Big) \log N+\frac{ N^{2} }{4}\Big\} .
\end{align*} 

\begin{lem}
For $l \geq 1 $, we have
\begin{align}
\prod_{ r=0 }^{ l-1 } \Gamma _{2} \big( \frac{ z+r }{ l } \big) ^{2l} \Gamma \big( \frac{ z+r }{ l } \big) ^{ 2r-2l+1 }
= k_{2} (l)^{-2} k_{1} (l) \cdot ( 2 \pi )^{ -lz } \cdot l^{ -z + z^{2} } \cdot \frac{ ( 2 \pi )^{z} \Gamma_{2} (z) ^{2} }{ \Gamma (z) } .
\label{Lem332A}
\end{align}
\end{lem}

\begin{proof}
We see that 
\begin{align*}
\prod_{ r_{1} = 0 }^{ l-1 } \prod_{ r_{2} = 0 }^{ l-1 } \Gamma_{2} \big( \frac{ z + r_{1} + r_{2} }{l} \big) 
& = \prod_{ r = 0 }^{ l-1 } \Gamma_{2} \big( \frac{ z + r}{l} \big) ^{ r+1 } \Gamma_{2} \big( \frac{ z + r }{l} \big) ^{ -r+l-1 } \\
& = \prod_{ r = 0 }^{ l-1 } \Gamma_{2} \big( \frac{ z + r}{l} \big) ^{l} \Gamma \big( \frac{ z + r }{l} \big) ^{ r-l+1 } .
\end{align*}
It follows from (\ref{GaMult}) and (\ref{DGMult}) that 
\begin{align*}
\prod_{ r=0 }^{ l-1 } \Gamma _{2} \big( \frac{ z+r }{ l } \big) ^{2l} \Gamma \big( \frac{ z+r }{ l } \big) ^{ 2r-2l+1 }
& = \prod_{ r_{1} = 0 }^{ l-1 } \prod_{ r_{2} = 0 }^{ l-1 } \Gamma_{ 2 } \big( \frac{ z + r_{1} + r_{2} }{2} \big) ^{2} \cdot
\prod_{ r=0 }^{ l-1 } \Gamma \big( \frac{ z+r }{l} \big) ^{ -1 } \\
& = k_{2} (l)^{-2} k_{1} (l) \cdot ( 2 \pi )^{ -lz } \cdot l^{ -z + z^{2} } \cdot 
\frac{ ( 2 \pi )^{z} \Gamma_{2} (z) ^{2} }{ \Gamma (z)} .
\end{align*}
This completes the proof. 
\end{proof}

From this lemma, we see that
\begin{align}
& \frac{d}{ds} \log \prod_{ r=0 }^{ l-1 } 
\Gamma _{2} \big( \frac{ s + \frac{l}{2} +r }{ l } \big) ^{2l} 
\Gamma \big( \frac{ s+ \frac{l}{2} +r }{ l } \big) ^{ 2r-2l+1 } \label{LDDGG} \\
= & - l \log ( 2 \pi ) + (-1+2s+l) \log l 
+ \frac{d}{ ds } \log \frac{ ( 2 \pi )^{s+ \frac{l}{2}} \Gamma_{2} ( s+ \frac{l}{2} ) ^{2} }{ \Gamma (s+ \frac{l}{2}) } . \notag
\end{align}
For simplicity, we put 
\begin{align*}
F_{l} (s) := & \exp \Big\{ 
- \frac{s}{l} \log \Big( k_{2} (l) ^{-2}  k_{1} (l) \Big) + 
\frac{ s(s+l) }{2} \log ( 2 \pi ) + \frac{ s(s+l)(3-2s-l) }{ 6l } \log l \\
& \qquad \; \; + \log \; \prod_{ r=0 }^{ l-1 } \; 
\Gamma _{3} \big( \frac{ s+l+r }{ l } \big) ^{-2l} \Gamma_{2} \big( \frac{ s+l+r }{ l } \big) ^{ -2r+2l-1 } \Big\}. \notag
\end{align*}
Then it is shown that $ F_{l} (s) $ satisfies
\begin{align}
\frac{ F_{l} (s)}{ F_{l} (s-l) } 
= & \Big( k_{2} (l) ^{-2}  k_{1} (l) \Big) ^{-1} \cdot ( 2 \pi )^{ls} \cdot l^{ s - s^{2} } \cdot \prod_{ r=0 }^{ l-1 } 
\Gamma _{2} \big( \frac{ s+r }{l} \big) ^{2l} \Gamma \big( \frac{ s+r }{ l } \big) ^{ 2r-2l+1 } \label{FSSML}  \\
= & \frac { ( 2 \pi )^{s} \Gamma_{2} (s) ^{2} }{ \Gamma (s) } , \notag
\end{align}
by (\ref{Lem332A}), and 
\begin{align}
\frac{ d^{2} }{ ds^{2} } \log F_{l} (s) 
= & \log ( 2 \pi ) + ( -2s +1-l ) \frac{ \log l }{ l } \label{LDDF} \\
& + \frac{ d^{2} }{ ds^{2} } \log \prod_{ r=0 }^{ l-1 }  
\Gamma _{3} \big( \frac{ s+l+r }{ l } \big) ^{-2l} \Gamma_{2} \big( \frac{ s+l+r }{ l } \big) ^{ -2r+2l-1 } . \notag
\end{align}

Combining (\ref{LDDGG}) and (\ref{LDDF}), we conclude that
\begin{align}
& \frac{ d^{2} }{ ds^{2} } \log \prod_{ r=0 }^{ l-1 } 
\Gamma _{3} \big( \frac{s+l+r}{l} \big) ^{ -2l } \Gamma _{2} \big( \frac{s+l+r}{l} \big) ^{ -2r + 2l -1 } \label{GaToF}\\
& \quad -  \frac{1}{l} \frac{d}{ ds } \log \prod_{ r=0 }^{ l-1 }
\Gamma _{2} \big( \frac{s+ \frac{l}{2} +r}{l} \big) ^{ 2l } \Gamma \big( \frac{s+ \frac{l}{2} +r}{l} \big) ^{ 2r - 2l +1 } \notag \\
=&  \frac{ d^{2} }{ ds^{2} } \log F_{l} (s) 
- \frac{1}{l} \frac{d}{ ds } \log \frac{ ( 2 \pi )^{s+ \frac{l}{2}} \Gamma_{2} ( s+ \frac{l}{2} ) ^{2} }{ \Gamma (s+ \frac{l}{2}) } . \notag
\end{align}

Now we have the description following. 

\begin{prop}
If $ h(r) = \phi _{s} (r)  + \phi _{s} ( -r ) $, we have
\begin{align}
\frac{ \mathrm{vol} ( \varGamma \setminus H ) }{ 4 \pi } \int_{ - \infty } ^{ \infty } h(r) r \tanh ( \pi r) dr 
= \frac{ d^{2} }{ ds^{2} } \log \Xi _{ l \infty , \mathrm{I} } (s) 
- \frac{1}{l} \frac{d}{ ds } \log \Xi _{ \mathrm{I} } ( s + \frac{l}{2} ) . \label{IC}
\end{align}
Here $ \Xi _{ l \infty , \mathrm{I} } (s) $ is defined by
\begin{align*}
\Xi _{ l \infty , \mathrm{I} } (s) := \exp \Big\{ \frac{ \mathrm{vol} ( \varGamma \setminus H ) }{ 2 \pi } \log F_{l} (s) \Big\} ,
\end{align*}
which satisfies
\begin{align}
\Xi _{ l \infty , \mathrm{I} } (s) = \frac { \Xi _{ l \infty , \mathrm{I} } (s+l) }{ \Xi _{ \mathrm{I} } (s+l) } , \quad 
\mathrm{Re} (s) > 1-l . \label{HIO}
\end{align} 
\end{prop}

\begin{proof} 
The first equation (\ref{IC}) follows immediately from (\ref{SelZetaGammaFactor}) and (\ref{GaToF}). 
The second one (\ref{HIO}) is clear from (\ref{FSSML}).
\end{proof}

\subsection{Elliptic Factor}

When $ \varGamma $ is the congruence subgroups of $ SL(2,\mathbb{Z}) $, 
the contribution of the elliptic and parabolic element in $ \varGamma $ appears 
in the trace formula. 
We now determine the elliptic factor. 
To begin with, we calculate the contribution of the elliptic conjugacy classes with order $2$. 

\begin{lem}
If $ h(r) = \phi _{s} (r)  + \phi _{s} ( -r ) $, we have
\begin{align} 
\frac{1}{2} \int_{ -\infty }^{ \infty } h(r) \frac{1}{ e^{ \pi r } + e^{ - \pi r} } dr 
= \frac{ d^{2} }{ ds^{2} } \log G_{l} (s) 
- \frac{1}{l} \frac{d}{ds} \log \Gamma \big( \frac{ s+ \frac{l}{2} }{2} \big) ^{-1} \Gamma \big( \frac{ s+ \frac{l}{2}+1 }{2}\big) ,
\label{Ell2}
\end{align}
where $ G_{l} (s) $ is given by
\begin{align*}
G_{l} (s) := 
\begin{cases} 
\displaystyle 
\exp \big( \frac{s}{2l} \log l \big) \cdot \prod_{ r=0 }^{ l-1 } \Gamma \big( \frac{ s+l+r }{ 2l } \big) ^{ (-1) ^{r} } 
& l: \text{odd}, \\ 
\displaystyle 
\exp \big( \frac{s}{2l} \log \frac{l}{2} \big) \cdot \prod_{ r=0 }^{ l-1 } \Gamma_{2} \big( \frac{ s+l+r }{l} \big) ^{ (-1) ^{r} } 
& l: \text{even} . 
\end{cases} 
\end{align*}
\end{lem}

\begin{proof}
Using the residue theorem, we see that
\begin{align*}
& \frac{1}{2} \int_{ -\infty }^{ \infty } h(r) \frac{1}{ e^{ \pi r } + e^{ - \pi r} } dr 
= \int_{ -\infty }^{ \infty } \phi _{s} (r) \frac{1}{ e^{ \pi r } + e^{ - \pi r} } dr 
= \sum_{ n=0 }^{ \infty } (-1)^{n} \phi_{s} \big( - \frac{ 2n+1 }{2} i \big) \\
= & \sum_{ n=0 }^{ \infty } \sum_{ m=0 }^{ \infty } (-1)^{n} 
\Big\{ \frac{1}{ (s+lm+l+n )^{2} } - \frac{1}{l} \frac{1}{ s+lm + \frac{l}{2} +n } + \frac{1}{l} \frac{1}{ s+lm + \frac{3l}{2} +n } \Big\} \\
= & \sum_{ r=0 }^{ l-1 } \sum_{ j=0 }^{ \infty } \sum_{ m=0 }^{ \infty } ( -1 )^{lj+r} \\
& \qquad \times \Big\{ \frac{1}{ (s+lm+l+lj+r )^{2} } - \frac{1}{l} \frac{1}{ s+lm+ \frac{l}{2} +lj+r } 
+ \frac{1}{l} \frac{1}{ s+lm + \frac{ 3l }{2} +lj+r } \Big\} \\
= & \sum_{ r=0 }^{ l-1 } (-1) ^{r} \sum_{ k=0 }^{ \infty } \sum_{ j=0 }^{k} (-1) ^{lj} 
\Big\{ \frac{1}{ (s+lk+l+r )^{2} } - \frac{1}{l} \frac{1}{ s+lk+ \frac{l}{2} +r } 
+ \frac{1}{l} \frac{1}{ s+lk + \frac{ 3l }{2} +r } \Big\} . \\
\end{align*}

Now we consider the following two cases. \\

\noindent \textit{ Case 1 :} $ l \equiv 1 \; ( \bmod \, 2 ) $.

If $l$ is odd, we observe that
\begin{align*}
\sum_{ j=0 }^{k} (-1) ^{lj} = 
\begin{cases} 
1 & k : \text{even}, \\ 
0 & k : \text{odd}. 
\end{cases} 
\end{align*}
Then we have
\begin{align*}
& \sum_{ r=0 }^{ l-1 } (-1) ^{r} \sum_{ k=0 }^{ \infty } \sum_{ j=0 }^{k} (-1) ^{lj} 
\Big\{ \frac{1}{ (s+lk+l+r )^{2} } - \frac{1}{l} \frac{1}{ s+lk+ \frac{l}{2} +r } 
+ \frac{1}{l} \frac{1}{ s+lk + \frac{ 3l }{2} +r } \Big\} \\
= &  \sum_{ r=0 }^{ l-1 } (-1) ^{r} \sum_{ k=0 }^{ \infty }  
\Big\{ \frac{1}{ (s+2lk+l+r )^{2} } - \frac{1}{l} \frac{1}{ s+2lk+ \frac{l}{2} +r } 
+ \frac{1}{l} \frac{1}{ s+2lk + \frac{ 3l }{2} +r } \Big\} . 
\end{align*}
Here, from the Weierstrass product expression for $ \Gamma (z) $, we obtain 
\begin{align*}
& \frac{ d^{2} }{ ds^{2} } \log \Gamma \big( \frac{ s+l+r }{ 2l } \big)
= \sum_{ k=0 }^{ \infty } \frac{1}{ (s+2lk+l+r) ^{2} } , \\
& \frac{d}{ ds } \log \Gamma \big( \frac{ s + \frac{l}{2} + r }{ 2l } \big) ^{-1} \Gamma \big( \frac{ s+ \frac{3l}{2} +r }{2l} \big) 
= \sum_{ k=0}^{ \infty } \Big\{ \frac{1}{ s+ 2lk+ \frac{l}{2} +r } - \frac{1}{ s+2lk+ \frac{3l}{2} +r } \Big\} . 
\end{align*}
This yields 
\begin{align*}
&  \sum_{ r=0 }^{ l-1 } (-1) ^{r} \sum_{ k=0 }^{ \infty }  
\Big\{ \frac{1}{ (s+2lk+l+r )^{2} } - \frac{1}{l} \frac{1}{ s+2lk+ \frac{l}{2} +r } 
+ \frac{1}{l} \frac{1}{ s+2lk + \frac{ 3l }{2} +r } \Big\} \\
= & \frac{ d^{2} }{ ds^{2} } \log \Big\{ \prod_{ r=0 }^{ l-1 } \Gamma \big( \frac{ s+l+r }{ 2l } \big) ^{(-1) ^{r}} \Big\} 
- \frac{1}{l} \frac{d}{ ds } \log \Big\{ \prod_{ r=0 }^{ l-1 } 
\Gamma \big( \frac{ s+ \frac{l}{2} +r }{ 2l } \big) ^{ -(-1)^{r} } \Gamma \big( \frac{ s+ \frac{3l}{2} +r }{2l} \big) ^{ (-1)^{r} } \Big\} .
\end{align*}
On the other hand, it is seen from (\ref{GaMult}) that 
\begin{align*}
\prod_{ r=0 } ^{ l-1 } \Gamma \big( \frac{ z+r }{ 2l } \big) ^{ - (-1)^{r} } \Gamma \big( \frac{ z+l+r }{ 2l } \big) ^{ (-1)^{r} } 
= & \prod_{ j=0 }^{ l-1 } \Gamma \big( \frac{ z+2j }{ 2l } \big) ^{ -1 } \Gamma \big( \frac{ z+1+ 2j }{ 2l } \big) \\
= & l^{ - \frac{1}{2} } \Gamma \big( \frac{z}{2} \big) ^{ -1 } \Gamma \big( \frac{z+1}{2} \big) . \notag
\end{align*}
Hence we arrive at
\begin{align*}
& \frac{ d^{2} }{ ds^{2} } \log \Big\{ \prod_{ r=0 }^{ l-1 } \Gamma \big( \frac{ s+l+r }{ 2l } \big) ^{(-1) ^{r}} \Big\} 
- \frac{1}{l} \frac{d}{ ds } \log \Big\{ \prod_{ r=0 }^{ \infty } 
\Gamma \big( \frac{ s + \frac{l}{2} + r }{ 2l } \big) ^{ - (-1)^{r} } \Gamma \big( \frac{ s+ \frac{3l}{2} +r }{2l} \big) ^{ (-1)^{r} } 
\Big\} \\
= & \frac{ d^{2} }{ ds^{2} } \log 
\Big\{ \exp \big( \frac{s}{2l} \log l \big) \prod_{ r=0 }^{ l-1 } \Gamma \big( \frac{ s+l+r }{ 2l } \big) ^{(-1) ^{r}} \Big\} 
- \frac{1}{l} \frac{d}{ ds } \log 
\Gamma \big( \frac{ s + \frac{l}{2} }{ 2 } \big) ^{ -1 } \Gamma \big( \frac{ s+ \frac{3l}{2} }{2} \big) . 
\end{align*} \\

\noindent \textit{ Case 2 :} $ l \equiv 2 \; ( \bmod \, 2 ) $.

Since we have 
\begin{align*}
\sum _{ j=0 }^{k} (-1) ^{ lj } = k+1 ,
\end{align*}
it follows that 
\begin{align*}
& \sum_{ r=0 }^{ l-1 } (-1) ^{r} \sum_{ k=0 }^{ \infty } \sum_{ j=0 }^{k} (-1) ^{lj} 
\Big\{ \frac{1}{ (s+lk+l+r )^{2} } - \frac{1}{l} \frac{1}{ s+lk+ \frac{l}{2} +r } 
+ \frac{1}{l} \frac{1}{ s+lk + \frac{ 3l }{2} +r } \Big\} \\
= &  \sum_{ r=0 }^{ l-1 } (-1) ^{r} \sum_{ k=0 }^{ \infty }  
\Big\{ \frac{ k+1 }{ (s+lk+l+r )^{2} } - \frac{1}{l} \frac{ k+1 }{ s+lk+ \frac{l}{2} +r } 
+ \frac{1}{l} \frac{ k+1 }{ s+lk + \frac{ 3l }{2} +r } \Big\} . 
\end{align*}
Here, from the Weierstrass product expression for $ \Gamma_{2} (z) $ and $ \Gamma (z) $ , we obtain 
\begin{align*}
& \frac{ d^{2} }{ ds^{2} } \log \Gamma_{2} \big( \frac{ s+l+r }{ l } \big)
= \frac{1}{ l^{2} } ( 1 + \gamma ) 
+ \sum_{ k=0 }^{ \infty } \big\{ \frac{ k+1 }{ (s+lk+l+r) ^{2} } - \frac{1}{ l^{2} } \frac{1}{k+1} \big\}, \\
& \frac{d}{ ds } \log \Gamma \big( \frac{ s + \frac{l}{2} + r }{ l } \big)  
= - \frac{ \gamma }{l}
- \sum_{ k=0}^{ \infty } \Big\{ \frac{1}{ s+ lk+ \frac{l}{2} +r } - \frac{1}{l} \frac{1}{k+1} \Big\} .
\end{align*}
This shows
\begin{align}
& \sum_{ k=0 }^{ \infty } \Big\{ \frac{ k+1 }{ (s+lk+l+r )^{2} } - \frac{1}{l} \frac{ k+1 }{ s+lk+ \frac{l}{2} +r } 
+ \frac{1}{l} \frac{k+1}{ s+lk + \frac{ 3l }{2} +r } \Big\} \label{DDGG} \\
& - \frac{ d^{2} }{ ds^{2} } \log \Gamma_{2} \big( \frac{ s+l+r }{ l } \big) 
- \frac{1}{l} \frac{d}{ ds } \log \Gamma \big( \frac{ s + \frac{l}{2} + r }{ l } \big) \notag \\
= & - \frac{1}{ l^{2} } 
+ \sum_{ k=0 }^{ \infty } \big\{ - \frac{1}{l} \frac{k}{ s+lk+ \frac{l}{2} +r } + \frac{1}{l} \frac{ k+1 }{ s+lk+ \frac{3l}{2} +r } 
\big\} \notag \\
= & - \frac{1}{ l^{2} } + \frac{1}{ l^{2} } ( s+ \frac{l}{2} +r ) 
\sum_{ k=0 }^{ \infty } \big\{ \frac{1}{ s+lk+ \frac{l}{2} +r } - \frac{1}{ s+lk+ \frac{3l}{2} +r } \big\} = 0 . \notag
\end{align}
Thus we see that 
\begin{align*}
& \sum_{ r=0 }^{ l-1 } (-1) ^{r} \sum_{ k=0 }^{ \infty }  
\Big\{ \frac{ k+1 }{ (s+lk+l+r )^{2} } - \frac{1}{l} \frac{ k+1 }{ s+lk+ \frac{l}{2} +r } 
+ \frac{1}{l} \frac{ k+1 }{ s+lk + \frac{ 3l }{2} +r } \Big\} \\
= & \frac{ d^{2} }{ ds^{2} } \log \prod_{ r=0 }^{ l-1 } \Gamma_{2} \big( \frac{ s+l+r }{ l } \big) ^{ (-1)^{r} } 
- \frac{1}{l} \frac{d}{ ds } \log \prod_{ r=0 }^{ l-1 } \Gamma \big( \frac{ s + \frac{l}{2} + r }{ l } \big) ^{ -(-1)^{r} } . \\
\end{align*}
Furthermore, since we have 
\begin{align*}
\prod_{ r=0 } ^{ l-1 } \Gamma \big( \frac{ z+r }{l} \big) ^{ - (-1)^{r} }  
= & \prod_{ j=0 }^{ \frac{l}{2} -1 } \Gamma \big( \frac{ z+2j }{ l } \big) ^{ -1 } \Gamma \big( \frac{ z+1+ 2j }{l} \big) \\
= & \big( \frac{l}{2} \big) ^{ - \frac{1}{2} } \Gamma \big( \frac{z}{2} \big) ^{ -1 } \Gamma \big( \frac{z+1}{2} \big) , 
\quad \text{ by (\ref{GaMult}), }  \notag
\end{align*}
we conclude that
\begin{align*}
& \frac{ d^{2} }{ ds^{2} } \log \prod_{ r=0 }^{ l-1 } \Gamma_{2} \big( \frac{ s+l+r }{ l } \big) ^{ (-1)^{r} } 
- \frac{1}{l} \frac{d}{ ds } \log \prod_{ r=0 }^{ l-1 } \Gamma \big( \frac{ s + \frac{l}{2} + r }{ l } \big) ^{ -(-1)^{r} }  \\
= & \frac{ d^{2} }{ ds^{2} } \log 
\Big\{ \exp \big( \frac{s}{2l} \log \frac{l}{2} \big) \prod_{ r=0 }^{ l-1 } \Gamma_{2} \big( \frac{ s+l+r }{ l } \big) ^{(-1) ^{r}} \Big\} 
- \frac{1}{l} \frac{d}{ ds } \log 
\Gamma \big( \frac{ s + \frac{l}{2} }{ 2 } \big) ^{ -1 } \Gamma \big( \frac{ s+ \frac{l}{2} +1 }{2} \big) . 
\end{align*}
This completes the proof. 
\end{proof}

We note that $ G_{l} (s) $ satisfies
\begin{align}
\frac{ G_{l} (s) }{ G_{l} (s-l) } = \Gamma \big( \frac{s}{2} \big) ^{ -1 } \Gamma \big( \frac{s+1}{2} \big) . \label{GlFu}
\end{align} 

Next, we calculate the contribution of elliptic conjugacy classes with order $3$.

\begin{lem}
If $ h(r) = \phi _{s} (r)  + \phi _{s} ( -r ) $, we have
\begin{align} 
\frac{1}{2 \sqrt{3} } \int_{ -\infty }^{ \infty } h(r) 
\frac{ e^{ \frac{\pi r}{3} } + e^{ - \frac{\pi r}{3} } }{ e^{ \pi r } + e^{ - \pi r} } dr 
= \frac{ d^{2} }{ ds^{2} } \log H_{l} (s) 
- \frac{1}{l} \frac{d}{ds} \log \Gamma \big( \frac{ s+ \frac{l}{2} }{3} \big) ^{-1} \Gamma \big( \frac{ s+ \frac{l}{2} +2 }{3} \big) ,
\label{Ell3}
\end{align}
where $ H_{l} (s) $ is given as follows. \\

For $ l \equiv 1 \; ( \bmod \, 3 ) $, we put 
\begin{align*}
H_{l} (s) := \exp \big( \frac{2s}{3l} \log l \big) 
\cdot & \prod_{ r=0 }^{ \frac{ l-1 }{3} } \Gamma \big( \frac{ s+l+3r }{ 3l } \big) \Gamma \big( \frac{ s+2l+3r }{ 3l } \big)  \\
\times & 
\prod_{ r=0 }^{ \frac{ l-4 }{3} } 
\Gamma \big( \frac{ s+l+3r+2 }{ 3l } \big) ^{-1} \Gamma \big( \frac{ s+2l+3r+1 }{ 3l } \big) ^{-1} , \notag
\end{align*}
for $ l \equiv 2 \; ( \bmod \, 3 ) $; 
\begin{align*}
H_{l} (s) := \exp \big( \frac{2s}{3l} \log l \big) 
\cdot & \prod_{ r=0 }^{ \frac{ l-2 }{3} } \Gamma \big( \frac{ s+l+3r }{ 3l } \big) \Gamma \big( \frac{ s+2l+3r+1 }{ 3l } \big)  \\
\times & 
\prod_{ r=0 }^{ \frac{ l-5 }{3} } \Gamma \big( \frac{ s+l+3r+2 }{ 3l } \big) ^{-1} \Gamma \big( \frac{ s+2l+3r+2 }{ 3l } \big) ^{-1} ,
\end{align*}
and for $ l \equiv 3 \; ( \bmod \, 3 ) $; 
\begin{align*}
H_{l} (s) := \exp \big( \frac{2s}{3l} \log \frac{l}{3} \big) 
\cdot \prod_{ r=0 }^{ \frac{ l-3 }{3} } \Gamma_{2} \big( \frac{ s+l+3r }{l} \big) \Gamma_{2} \big( \frac{ s+l+3r+2 }{l} \big)^{-1} .
\end{align*}
\end{lem}

\begin{proof}
Applying the residue theorem as usual, we calculate as
\begin{align*}
& \frac{1}{ 2\sqrt{3} } \int_{ -\infty }^{ \infty } h(r) 
\frac{ e^{ \frac{\pi r}{3} } + e^{ - \frac{\pi r}{3} } }{ e^{ \pi r } + e^{ - \pi r} } dr 
= \frac{1}{ \sqrt{3} } \sum_{ n=0 }^{ \infty } 
2 (-1)^{n} \cos \big\{ \frac{ \pi }{3} \big( n + \frac{1}{2} \big) \big\} \cdot \phi_{s} \big( - \frac{ 2n+1 }{2} i \big) \\
\\
= & \frac{2}{ \sqrt{3} } \sum_{ r=0 }^{ l-1 } \sum_{ k=0 }^{ \infty } \sum_{ j=0 }^{k} 
(-1)^{ lj+r } \cos \big\{ \frac{ \pi }{3} \big( lj+r + \frac{1}{2} \big) \big\} \\
& \qquad \times \Big\{ \frac{1}{ (s+lk+l+r )^{2} } - \frac{1}{l} \frac{1}{ s+lk+ \frac{l}{2} +r } 
+ \frac{1}{l} \frac{1}{ s+lk+ \frac{ 3l }{2} +r } \Big\} .
\end{align*}

Now we consider the following three cases. \\

\noindent \textit{ Case 1 :} $ l \equiv 1 \; ( \bmod \, 3 ) $.

If $ l \equiv 1 \; ( \bmod \, 3 ) $, we have 
\begin{align*}
\sum_{ j=0 }^{k} (-1)^{ lj+r } \cos \big\{ \frac{ \pi }{3} \big( lj+r + \frac{1}{2} \big) \big\} 
= \begin{cases}
(-1)^{r} \cos \big\{ \frac{ \pi }{3} \big( r+ \frac{1}{2} \big) \big\} & k \equiv 0 \; ( \bmod \, 3 ) , \\
(-1)^{r} \cos \big\{ \frac{ \pi }{3} \big( r- \frac{1}{2} \big) \big\} & k \equiv 1 \; ( \bmod \, 3 ) , \\
0 & k \equiv 2 \; ( \bmod \, 3 ) .
\end{cases}
\end{align*}
Then we see that
\begin{align*}
& \frac{2}{ \sqrt{3} } \sum_{ r=0 }^{ l-1 } \sum_{ k=0 }^{ \infty } \sum_{ j=0 }^{k} 
(-1)^{ lj+r } \cos \big\{ \frac{ \pi }{3} \big( lj+r + \frac{1}{2} \big) \big\} \\
& \qquad \times \Big\{ \frac{1}{ (s+lk+l+r )^{2} } - \frac{1}{l} \frac{1}{ s+lk+ \frac{l}{2} +r } 
+ \frac{1}{l} \frac{1}{ s+lk+ \frac{ 3l }{2} +r } \Big\} \\
= & \frac{2}{ \sqrt{3} } \sum_{ r=0 }^{ l-1 } (-1)^{r} \cos \big\{ \frac{ \pi }{3} \big( r+ \frac{1}{2} \big) \big\} \\
& \qquad \times \sum_{ k=0 }^{ \infty } \Big\{ \frac{1}{ (s+3lk+l+r )^{2} } - \frac{1}{l} \frac{1}{ s+3lk+ \frac{l}{2} +r } 
+ \frac{1}{l} \frac{1}{ s+3lk+ \frac{ 3l }{2} +r } \Big\} \\
+ & \frac{2}{ \sqrt{3} } \sum_{ r=0 }^{ l-1 } (-1)^{r} \cos \big\{ \frac{ \pi }{3} \big( r - \frac{1}{2} \big) \big\} \\
& \qquad \times \sum_{ k=0 }^{ \infty } \Big\{ \frac{1}{ (s+3lk+2l+r )^{2} } - \frac{1}{l} \frac{1}{ s+3lk+ \frac{3l}{2} +r } 
+ \frac{1}{l} \frac{1}{ s+3lk+ \frac{5l}{2} +r } \Big\} .
\end{align*}
Here we observe that
\begin{align}
\frac{2}{ \sqrt{3} }(-1)^{r} \cos \big\{ \frac{ \pi }{3} \big( r + \frac{1}{2} \big) \big\} 
= \begin{cases}
1  & r \equiv 0 \; ( \bmod \, 3 ) , \\
0  & r \equiv 1 \; ( \bmod \, 3 ) , \\
-1 & r \equiv 2 \; ( \bmod \, 3 ) .
\end{cases} \label{CCSS}
\end{align}
In addition, since we have
\begin{align*}
& \frac{ d^{2} }{ ds^{2} } \log \Gamma \big( \frac{ s+l+r }{ 3l } \big)
= \sum_{ k=0 }^{ \infty } \frac{1}{ (s+3lk+l+r) ^{2} } , \\
& \frac{d}{ ds } \log \Gamma \big( \frac{ s + \frac{l}{2} + r }{ 3l } \big) ^{-1} \Gamma \big( \frac{ s+ \frac{3l}{2} +r }{3l} \big) 
= \sum_{ k=0}^{ \infty } \Big\{ \frac{1}{ s+ 3lk+ \frac{l}{2} +r } - \frac{1}{ s+3lk+ \frac{3l}{2} +r } \Big\} , 
\end{align*}
it follows from (\ref{GaMult}) that
\begin{align*}
& \frac{2}{ \sqrt{3} } \sum_{ r=0 }^{ l-1 } (-1)^{r} \cos \big\{ \frac{ \pi }{3} \big( r+ \frac{1}{2} \big) \big\} \\
& \qquad \times \sum_{ k=0 }^{ \infty } \Big\{ \frac{1}{ (s+3lk+l+r )^{2} } - \frac{1}{l} \frac{1}{ s+3lk+ \frac{l}{2} +r } 
+ \frac{1}{l} \frac{1}{ s+3lk+ \frac{ 3l }{2} +r } \Big\} \\
+ & \frac{2}{ \sqrt{3} } \sum_{ r=0 }^{ l-1 } (-1)^{r} \cos \big\{ \frac{ \pi }{3} \big( r - \frac{1}{2} \big) \big\} \\
& \qquad \times \sum_{ k=0 }^{ \infty } \Big\{ \frac{1}{ (s+3lk+2l+r )^{2} } - \frac{1}{l} \frac{1}{ s+3lk+ \frac{3l}{2} +r } 
+ \frac{1}{l} \frac{1}{ s+3lk+ \frac{5l}{2} +r } \Big\} \\
\\
= & \frac{ d^{2} }{ ds^{2} } \log 
\prod_{r=0}^{ \frac{l-1}{3} } \Gamma \big( \frac{ s+l+3r }{ 3l } \big) \Gamma \big( \frac{ s+2l+3r }{ 3l } \big) \cdot
\prod_{r=0}^{ \frac{l-4}{3} } \Gamma \big( \frac{ s+l+3r+2 }{ 3l } \big) ^{-1} \Gamma \big( \frac{ s+2l+3r+1 }{ 3l } \big)^{-1} \\
& \qquad - \frac{1}{l} \frac{d}{ds} \log
\prod_{ r=0 }^{ l-1 } \Gamma \big( \frac{ s+ \frac{l}{2} +3r }{ 3l } \big) ^{-1} \Gamma \big( \frac{ s+\frac{l}{2}+2+3r }{ 3l } \big) \\
= & \frac{ d^{2} }{ ds^{2} } \log H_{l} (s)
- \frac{1}{l} \frac{d}{ds} \log \Gamma \big( \frac{ s+ \frac{l}{2} }{3} \big) ^{-1} \Gamma \big( \frac{ s+ \frac{l}{2} +2 }{3} \big) .
\end{align*} 

\noindent \textit{ Case 2 :} $ l \equiv 2 \; ( \bmod \, 3 ) $.

If $ l \equiv 2 \; ( \bmod \, 3 ) $, we see that
\begin{align*}
\sum_{ j=0 }^{k} (-1)^{ lj+r } \cos \big\{ \frac{ \pi }{3} \big( lj+r + \frac{1}{2} \big) \big\} 
= \begin{cases}
(-1)^{r} \cos \big\{ \frac{ \pi }{3} \big( r+ \frac{1}{2} \big) \big\} & k \equiv 0 \; ( \bmod \, 3 ) , \\
(-1)^{r} \cos \big\{ \frac{ \pi }{3} \big( r+ \frac{3}{2} \big) \big\} & k \equiv 1 \; ( \bmod \, 3 ) , \\
0 & k \equiv 2 \; ( \bmod \, 3 ) .
\end{cases}
\end{align*}
Using (\ref{CCSS}), we have
\begin{align*}
& \frac{2}{ \sqrt{3} } \sum_{ r=0 }^{ l-1 } \sum_{ k=0 }^{ \infty } \sum_{ j=0 }^{k} 
(-1)^{ lj+r } \cos \big\{ \frac{ \pi }{3} \big( lj+r + \frac{1}{2} \big) \big\} \\
& \qquad \times \Big\{ \frac{1}{ (s+lk+l+r )^{2} } - \frac{1}{l} \frac{1}{ s+lk+ \frac{l}{2} +r } 
+ \frac{1}{l} \frac{1}{ s+lk+ \frac{ 3l }{2} +r } \Big\} \\
= & \sum_{ r=0 }^{ \frac{l-2}{3} } \Big\{ 
\frac{ d^{2} }{ ds^{2} } \log \Gamma \big( \frac{ s+2l+3r+1 }{ 3l } \big) - \frac{1}{l} 
\frac{d}{ds} \log \Gamma \big( \frac{ s+ \frac{3l}{2} +3r+1 }{3l} \big) ^{-1} 
\Gamma \big( \frac{ s+ \frac{5l}{2} +3r+1 }{3l} \big) \Big\} \\
- & \sum_{ r=0 }^{ \frac{l-5}{3} } \Big\{ 
\frac{ d^{2} }{ ds^{2} } \log \Gamma \big( \frac{ s+2l+3r+2 }{ 3l } \big) - \frac{1}{l} 
\frac{d}{ds} \log \Gamma \big( \frac{ s+ \frac{3l}{2} +3r+2 }{3l} \big) ^{-1} 
\Gamma \big( \frac{ s+ \frac{5l}{2} +3r+2 }{3l} \big) \Big\} \\
= & \frac{ d^{2} }{ ds^{2} } \log 
\prod_{r=0}^{ \frac{l-2}{3} } \Gamma \big( \frac{ s+l+3r }{ 3l } \big) \Gamma \big( \frac{ s+2l+3r+1 }{ 3l } \big) \cdot
\prod_{r=0}^{ \frac{l-5}{3} } \Gamma \big( \frac{ s+l+3r+2 }{ 3l } \big) ^{-1} \Gamma \big( \frac{ s+2l+3r+2 }{ 3l } \big)^{-1} \\
& \qquad - \frac{1}{l} \frac{d}{ds} \log
\prod_{ r=0 }^{ l-1 } \Gamma \big( \frac{ s+ \frac{l}{2} +3r }{ 3l } \big) ^{-1} \Gamma \big( \frac{ s+\frac{l}{2}+2+3r }{ 3l } \big) \\
= & \frac{ d^{2} }{ ds^{2} } \log H_{l} (s)
- \frac{1}{l} \frac{d}{ds} \log \Gamma \big( \frac{ s+ \frac{l}{2} }{3} \big) ^{-1} \Gamma \big( \frac{ s+ \frac{l}{2} +2 }{3} \big) .
\end{align*} \\

\noindent \textit{ Case 3 :} $ l \equiv 3 \; ( \bmod \, 3 ) $.

If $ l \equiv 2 \; ( \bmod \, 3 ) $, we observe that
\begin{align*}
\sum_{ j=0 }^{k} (-1)^{ lj+r } \cos \big\{ \frac{ \pi }{3} \big( lj+r + \frac{1}{2} \big) \big\} 
= (k+1) \cdot (-1)^{r} \cos \big\{ \frac{ \pi }{3} \big( r+ \frac{1}{2} \big) \big\} .
\end{align*}
Then we see that 
\begin{align*}
& \frac{2}{ \sqrt{3} } \sum_{ r=0 }^{ l-1 } \sum_{ k=0 }^{ \infty } \sum_{ j=0 }^{k} 
(-1)^{ lj+r } \cos \big\{ \frac{ \pi }{3} \big( lj+r + \frac{1}{2} \big) \big\} \\
& \qquad \times \Big\{ \frac{1}{ (s+lk+l+r )^{2} } - \frac{1}{l} \frac{1}{ s+lk+ \frac{l}{2} +r } 
+ \frac{1}{l} \frac{1}{ s+lk+ \frac{ 3l }{2} +r } \Big\} \\
= & \frac{2}{ \sqrt{3} } \sum_{ r=0 }^{ l-1 } 
(-1)^{r} \cos \big\{ \frac{ \pi }{3} \big( r + \frac{1}{2} \big) \big\} \\
& \qquad \times \sum_{ k=0 }^{ \infty } \Big\{ \frac{ k+1 }{ (s+lk+l+r )^{2} } - \frac{1}{l} \frac{ k+1 }{ s+lk+ \frac{l}{2} +r } 
+ \frac{1}{l} \frac{ k+1 }{ s+lk+ \frac{ 3l }{2} +r } \Big\} .
\end{align*}
It follows from (\ref{DDGG}) and (\ref{CCSS}) that
\begin{align*}
& \frac{2}{ \sqrt{3} } \sum_{ r=0 }^{ l-1 } 
(-1)^{r} \cos \big\{ \frac{ \pi }{3} \big( r + \frac{1}{2} \big) \big\} \\
& \qquad \times \sum_{ k=0 }^{ \infty } \Big\{ \frac{ k+1 }{ (s+lk+l+r )^{2} } - \frac{1}{l} \frac{ k+1 }{ s+lk+ \frac{l}{2} +r } 
+ \frac{1}{l} \frac{ k+1 }{ s+lk+ \frac{ 3l }{2} +r } \Big\} \\
= & \sum_{r=0}^{ \frac{l-3}{3} } \Big\{ \frac{ d^{2} }{ ds^{2} } \log \Gamma_{2} \big( \frac{ s+l+3r }{l} \big) 
- \frac{1}{l} \frac{d}{ds} \log \Gamma \big( \frac{ s+ \frac{l}{2} +3r }{l} \big) \Big\} \\
& \quad - \sum_{r=0}^{ \frac{l-3}{3} } \Big\{ \frac{ d^{2} }{ ds^{2} } \log \Gamma_{2} \big( \frac{ s+l+3r+2 }{l} \big) ^{-1} 
- \frac{1}{l} \frac{d}{ds} \log \Gamma \big( \frac{ s+ \frac{l}{2} +3r+2 }{l} \big) ^{-1} \Big\} \\
= & \frac{ d^{2} }{ ds^{2} } \log 
\prod_{r=0}^{ \frac{l-3}{3} } \Gamma_{2} \big( \frac{ s+l+3r }{l} \big) \Gamma_{2} \big( \frac{ s+l+3r+2 }{l} \big) ^{-1} \\
& \qquad - \frac{1}{l} \frac{d}{ds} \log \prod_{ r=0 }^{ \frac{l-3}{3} } 
\Gamma \big( \frac{ s+ \frac{l}{2} +3r }{l} \big) ^{-1} \Gamma \big( \frac{ s+\frac{l}{2}+2+3r }{l} \big) .
\end{align*} 
Since we have
\begin{align*}
\prod_{ r=0 }^{ \frac{l-3}{3} } \Gamma \big( \frac{ z+3r }{l} \big) ^{-1} \Gamma \big( \frac{z+2+3r }{l} \big)
= \big( \frac{l}{3} \big) ^{ - \frac{2}{3} } \Gamma \big( \frac{z}{3} \big) ^{-1} \Gamma \big( \frac{ z+2 }{3} \big), 
\quad \text{by (\ref{GaMult}),}
\end{align*}
we conclude that
\begin{align*}
& \frac{ d^{2} }{ ds^{2} } \log 
\prod_{r=0}^{ \frac{l-3}{3} } \Gamma_{2} \big( \frac{ s+l+3r }{l} \big) \Gamma_{2} \big( \frac{ s+l+3r+2 }{l} \big) ^{-1} \\
& \qquad - \frac{1}{l} \frac{d}{ds} \log \prod_{ r=0 }^{ \frac{l-3}{3} } 
\Gamma \big( \frac{ s+ \frac{l}{2} +3r }{l} \big) ^{-1} \Gamma \big( \frac{ s+\frac{l}{2}+2+3r }{l} \big) \\
= & \frac{ d^{2} }{ ds^{2} } \log H_{l} (s)
- \frac{1}{l} \frac{d}{ds} \log \Gamma \big( \frac{ s+ \frac{l}{2} }{3} \big) ^{-1} \Gamma \big( \frac{ s+ \frac{l}{2} +2 }{3} \big) .
\end{align*} 
This completes the proof. 
\end{proof}

We see that $ H_{l} (s) $ satisfies
\begin{align}
\frac{ H_{l} (s) }{ H_{l} (s-l) } = \Gamma \big( \frac{s}{3} \big) ^{-1} \Gamma \big( \frac{s+2}{3} \big) . \label{HlFu}
\end{align}

Now we describe the contribution of the elliptic terms as follows. 

\begin{prop}
If $ h(r) = \phi _{s} (r)  + \phi _{s} ( -r ) $, we have
\begin{align}
& \frac{ \nu_{2} }{4} \int_{ - \infty }^{ \infty } h(r) \frac{1}{ e^{ \pi r } + e^{ - \pi r } }  dr 
+ \frac{ \nu_{3}  }{ 3 \sqrt{3} } \int_{ - \infty }^{ \infty } h(r) 
\frac{ e^{ \frac{ \pi r }{3} } + e^{ - \frac{ \pi r }{3} } }{ e^{ \pi r} + e^{ - \pi r } } dr \label{EC} \\
= & \frac{ d^{2} }{ ds^{2} } \log \Xi _{ l \infty , \mathrm{ell} } (s) 
- \frac{1}{l} \frac{d}{ ds } \log \Xi _{ \mathrm{ell} } ( s + \frac{l}{2} ) , \notag
\end{align}
where $ \Xi _{ l \infty , \mathrm{ell} } (s) $ is defined by
\begin{align*}
\Xi _{ l \infty , \mathrm{ell} } (s) 
:= \exp \big\{ \frac{ \nu_{2} }{2} \log G_{l} (s) + \frac{ 2 \nu_{3} }{3} \log H_{l} (s) \big\} .
\end{align*}
Furthermore, the function $ \Xi _{ l \infty , \mathrm{ell} } (s) $ satisfies 
\begin{align}
\Xi _{ l \infty , \mathrm{ell} } (s) = \frac { \Xi _{ l \infty , \mathrm{ell} } (s+l) }{ \Xi _{ \mathrm{ell} } (s+l) } , \quad 
\mathrm{Re} (s) > 1-l . \label{HEO}
\end{align} 
\end{prop}

\begin{proof}
This first equation follows immediately from (\ref{Ell2}) and (\ref{Ell3}). 
The second one is clear from (\ref{GlFu}) and (\ref{HlFu}). 
\end{proof}

\subsection{Parabolic Factor}

It remains to calculate the parabolic terms of the trace formula. 

\begin{lem}
If $ h(r) = \phi _{s} (r)  + \phi _{s} ( -r ) $, we have
\begin{align*}
- g(0) \log \big( \mathcal{A} \frac{ 2^{ \kappa } }{ \pi ^{\kappa} } \big) = 0 
= \frac{ d^{2} }{ ds^{2} } \log \Big( \mathcal{A} \frac{ 2^{ \kappa } }{ \pi ^{\kappa} } \Big) ^{ - \frac{ s(s+l) }{2l} }
- \frac{1}{l} \frac{d}{ds} \log \Big( \mathcal{A} \frac{ 2^{ \kappa } }{ \pi ^{\kappa} } \Big) ^{ - ( s+ \frac{l}{2} ) } ,
\end{align*}
and
\begin{align*}
\frac{1}{4} ( \kappa - \kappa _{0} ) h(0)
= \frac{ d^{2} }{ ds^{2} } \log \Big\{ \exp \big( \frac{s}{l} \log l \big) \cdot \Gamma \big( \frac{ s+l - \frac{1}{2} }{l} \big) \Big\}
^{\frac{ \kappa - \kappa _{0} }{2} } 
- \frac{1}{l} \frac{d}{ds} \log \big( s+ \frac{l}{2} - \frac{1}{2} \big) ^{\frac{ \kappa - \kappa _{0} }{2} } .
\end{align*}
\end{lem}

\begin{proof}
From (\ref{TestFun}), we see that
$ g(0) = 0 $, and
\begin{align*}
\frac{1}{2} h(0) =& \phi _{s} (0)
= \sum_{ m=0 }^{ \infty } \Big\{ 
\frac{1}{ ( s+lm+l - \frac{1}{2} )^{2} } \Big\} - \frac{1}{l} \frac{1}{ s+ \frac{l}{2} - \frac{1}{2} } \\
=& \frac{ d^{2} }{ ds^{2} } \log \Gamma \big( \frac{ s+l - \frac{1}{2} }{l} \big)  
- \frac{1}{l} \frac{d}{ds} \log \big( s+ \frac{l}{2} - \frac{1}{2} \big) \\
=& \frac{ d^{2} }{ ds^{2} } \log \Big\{ \exp \big( \frac{s}{l} \log l \big) 
\cdot \Gamma \big( \frac{ s+l - \frac{1}{2} }{l} \big) \Big\}  
- \frac{1}{l} \frac{d}{ds} \log \big( s+ \frac{l}{2} - \frac{1}{2} \big) .
\end{align*}
Hence, two equations of the lemma hold. 
\end{proof}

\begin{lem}
If $ h(r) = \phi _{s} (r)  + \phi _{s} ( -r ) $, we have
\begin{align}
& - \frac{1}{ 2 \pi } \int_{ - \infty }^{ \infty } h(r) \big\{ 
\frac{ \Gamma ^{ \prime } }{ \Gamma } ( 1+ir ) + \frac{ \Gamma ^{ \prime } }{ \Gamma } ( \frac{1}{2} +ir ) \big\} dr \\
= & \frac{ d^{2} }{ ds^{2} } \log \{ I_{l} ( s + \frac{1}{2} ) \cdot I_{l} (s) \}
- \frac{1}{l} \frac{d}{ds} \log \big\{
\Gamma \big( s+ \frac{l}{2} + \frac{1}{2} \big) \cdot \Gamma \big( s+ \frac{l}{2} \big) \big\} ^{-1} , \notag
\end{align}
where $ I_{l} (s) $ is defined by
\begin{align*}
I_{l} (s) & := \exp \big\{ - \frac{s}{l} \log \big( k_{1} (l) \big) - \frac{ s(s+l) }{2l} \log l \big\} \cdot
\prod _{ r=0 }^{ l-1 } \Gamma _{2} \big( \frac{ s+l+r }{l} \big) ,
\end{align*}
which satisfies $ I_{l} (s) / I_{l} (s-l)  = \Gamma (s) ^{-1} $.
\end{lem}

\begin{proof}
Applying the residue theorem to the lower half plane, we see that
\begin{align*}
& - \frac{1}{ 2 \pi } \int_{ - \infty }^{ \infty } h(r) 
\frac{ \Gamma ^{ \prime } }{ \Gamma } ( 1+ir ) dr 
= - \frac{1}{ 2 \pi } \int_{ - \infty }^{ \infty } \phi _{s} (r) 
\big\{ \frac{ \Gamma ^{ \prime } }{ \Gamma } ( 1+ir ) + \frac{ \Gamma ^{ \prime } }{ \Gamma } ( 1-ir ) \big\} dr \\
= & \sum_{ n=0 }^{ \infty } \phi _{s} \big( -i(n+1) \big) \\
= & \sum_{ n=0 }^{ \infty } \sum_{ m=0 }^{ \infty } 
\Big\{ \frac{1}{ ( s+lm+l+n + \frac{1}{2}  )^{2} } - \frac{1}{l} \frac{1}{ s+lm + \frac{l}{2} + n + \frac{1}{2} } 
+ \frac{1}{l} \frac{1}{ s+lm + \frac{3l}{2} +n + \frac{1}{2} } \Big\} \\
= & \sum_{ r=0 }^{ l-1 }  \sum_{ k=0 }^{ \infty }  
\Big\{ \frac{ k+1 }{ ( s+lk+l+r + \frac{1}{2}  )^{2} } - \frac{1}{l} \frac{ k+1 }{ s+lk+ \frac{l}{2} + r + \frac{1}{2} } 
+ \frac{1}{l} \frac{ k+1 }{ s+lk + \frac{3l}{2} + r + \frac{1}{2} } \Big\} .
\end{align*}
From (\ref{DDGG}), we have
\begin{align*}
& \sum_{ r=0 }^{ l-1 }  \sum_{ k=0 }^{ \infty }  
\Big\{ \frac{ k+1 }{ ( s+lk+l+r + \frac{1}{2}  )^{2} } - \frac{1}{l} \frac{ k+1 }{ s+lk+ \frac{l}{2} + r + \frac{1}{2} } 
+ \frac{1}{l} \frac{ k+1 }{ s+lk + \frac{3l}{2} + r + \frac{1}{2} } \Big\} \\
= & \sum_{ r=0 }^{ \infty } \Big\{ \frac{ d^{2} }{ ds^{2} } \log \Gamma_{2} \big( \frac{ s+l+r + \frac{1}{2} }{l} \big) 
- \frac{1}{l} \frac{d}{ds} \log \Gamma \big( \frac{ s + \frac{l}{2} +r + \frac{1}{2} }{l} \big)^{-1} \Big\} .
\end{align*}
It follows from the multiplication formula $( \ref{GaMult} )$ that
\begin{align*}
& \sum_{ r=0 }^{ \infty } \Big\{ \frac{ d^{2} }{ ds^{2} } \log \Gamma_{2} \big( \frac{ s+l+r + \frac{1}{2} }{l} \big) 
- \frac{1}{l} \frac{d}{ds} \log \Gamma \big( \frac{ s + \frac{l}{2} +r + \frac{1}{2} }{l} \big)^{-1} \Big\} \\
= & \frac{ d^{2} }{ ds^{2} } \log \prod_{ r=0 }^{ l-1 } \Gamma_{2} \big( \frac{ s+l+r + \frac{1}{2} }{l} \big)
- \frac{1}{l} \frac{d}{ds} \log 
\Big\{ k_{1} (l) \cdot l^{ s + \frac{l}{2} + \frac{1}{2} }\Gamma \big( s + \frac{l}{2} + \frac{1}{2} \big)^{-1} \Big\} \\
= & \frac{ d^{2} }{ ds^{2} } \log I_{l} ( s + \frac{1}{2} ) 
- \frac{1}{l} \frac{d}{ds} \log \Gamma \big( s + \frac{l}{2} + \frac{1}{2} \big)^{-1} .
\end{align*} 
Similarly we calculate as
\begin{align*}
& - \frac{1}{ 2 \pi } \int_{ - \infty }^{ \infty } h(r) 
\frac{ \Gamma ^{ \prime } }{ \Gamma } ( \frac{1}{2} +ir ) dr 
= \sum_{ n=0 }^{ \infty } \phi _{s} \big( -i( n+ \frac{1}{2} ) \big) \\
= & \sum_{ n=0 }^{ \infty } \sum_{ m=0 }^{ \infty } 
\Big\{ \frac{1}{ ( s+lm+l+n )^{2} } - \frac{1}{l} \frac{1}{ s+lm + \frac{l}{2} + n } 
+ \frac{1}{l} \frac{1}{ s+lm + \frac{3l}{2} + n } \Big\} \\
= & \frac{ d^{2} }{ ds^{2} } \log I_{l} (s) - \frac{1}{l} \frac{d}{ds} \log \Gamma \big( s + \frac{l}{2} \big)^{-1} 
\end{align*}
This proves the lemma. 
\end{proof}

As mentioned in Section 2, 
it is seen that 
the Dirichlet $L$-function appears in the parabolic factor of the Selberg zeta function. 
Here we introduce the higher Dirichlet $L$-function as follows. 

\begin{lem}
If $ h(r) = \phi _{s} (r)  + \phi _{s} ( -r ) $, we have
\begin{align}
2 \sum_{ n=1 }^{ \infty } \frac{ \chi (n) \Lambda (n) }{n} g(2 \log n) 
= \frac{ d^{2} }{ ds^{2} } \log L_{ 2l, \infty} ( 2s , \chi )
- \frac{1}{l} \frac{d}{ds} \log L ( 2s+l , \chi )^{-1} .
\end{align}
Here $ L_{ l \infty} (s , \chi) $ is the higher Dirichlet $L$-function defined by 
\begin{align*}
L_{ l\infty} (s , \chi) := \prod_{m=1}^{\infty} L (s+lm, \chi ) 
= \prod_{m=1}^{\infty} \prod_{ p: \text{prime} } \big(\ 1 - \chi(p) p^{-s-lm} \big).
\end{align*} 
\end{lem}

\begin{proof}
Since 
\begin{align*}
\frac{d}{ds} \log L ( s , \chi ) 
= - \sum_{ n=1 }^{ \infty } { \chi (n) \Lambda (n) }{ n^{s} } ,
\end{align*}
it follows that 
\begin{align*}
\frac{d}{ds} \log L_{ l \infty } ( s ,\chi ) = \sum_{m=1}^{\infty} \frac{d}{ds} \log L ( s+lm, \chi )
= - \sum_{m=1}^{\infty} \sum_{n=2}^{\infty} \frac{ \chi(n) \Lambda (n)}{ n^{s+lm} } 
= - \sum_{n=2}^{\infty} \frac{ \chi(n) \Lambda (n)}{ n^{l} -1 } \frac{1}{ n^{s} }, 
\end{align*}
and
\begin{align*}
\frac{ d^{2} }{ d s^{2} } \log L_{ l \infty } ( s ,\chi ) 
= \sum_{n=2}^{\infty} \frac{ \chi(n) \Lambda (n)}{ n^{l} -1 } \frac{ \log n }{ n^{s} } .
\end{align*}
Since we have
\begin{align*}
g ( 2 \log n ) = \frac{ 2 \log n }{ n^{2l} - 1 } \frac{1}{ n^{2s-1} } - \frac{1}{l} \frac{1}{ n^{2s+l-1} }, 
\end{align*}
the lemma follows. 
\end{proof}

The contribution of the parabolic terms are now described as follows. 

\begin{prop}
If $ h(r) = \phi _{s} (r)  + \phi _{s} ( -r ) $, we have
\begin{align}
&- g(0) \log \big( \mathcal{A} \frac{ 2^{ \kappa } }{ \pi ^{ \kappa } } \big) 
+ \frac{1}{4} ( \kappa - \kappa _{0} ) h(0) \label{} \label{PC} \\
&- \frac{ \kappa }{ 2 \pi } \int_{ - \infty }^{ \infty } h(r) 
\big\{ \frac{ \Gamma ^{ \prime } }{ \Gamma }  (1 + ir) +  \frac{ \Gamma ^{ \prime } }{ \Gamma } ( \frac{1}{2}  + ir ) \big\} dr 
+ 2 \sum_{ \chi } \sum_{ n=1 }^{ \infty } \frac{ \chi (n) \Lambda (n) }{n} g(2 \log n) \notag \\
= & \frac{ d^{2} }{ ds^{2} } \log \Xi _{ l \infty , \mathrm{par} } (s) 
- \frac{1}{l} \frac{d}{ ds } \log \Xi _{ \mathrm{par} } ( s + \frac{l}{2} ) . \notag
\end{align}
Here $ \Xi _{ l \infty , \mathrm{par} } (s) $ is defined by
\begin{align*}
\Xi _{ l \infty , \mathrm{par} } (s) 
:= &  \exp \Big\{ - \frac{ s(s+l) }{2l} \log \Big( \mathcal{A} \frac{ 2^{ \kappa } }{ \pi ^{\kappa} } \Big) 
+ \frac{ \kappa - \kappa _{0} }{2} \Big( 
\frac{s}{l} \log l + \log \Gamma \big( \frac{ s+l - \frac{1}{2} }{l} \big) \Big) \Big\} \notag \\ 
\times & \Big\{ I_{l} ( s + \frac{1}{2} ) \cdot I_{l} (s) \Big\} ^{ \kappa } \cdot \prod_{ \chi } L_{ 2l, \infty } ( 2s , \chi) , \notag
\end{align*}
which satisfies 
\begin{align}
\Xi _{ l \infty , \mathrm{par} } (s) = \frac { \Xi _{ l \infty , \mathrm{par} } (s+l) }{ \Xi _{ \mathrm{par} } (s+l) } . \label{HPO}
\end{align} 
\qed
\end{prop}

\subsection{Complete Higher Selberg Zeta Function}

We have already calculated all terms of the trace formula except for the spectral terms. 
Using these calculations, we show the functional equation of higher Selberg zeta function for congruence subgroups. 
To describe the functional equation in symmetric way, we introduce a periodic function
\begin{align}
\Theta _{ l \infty } (s)
:= \prod_{ n=0 }^{ \infty }
\big( 1 - e^{ \frac{ 2 \pi i }{l} ( \frac{1}{2} + i r_{n} - s ) } \big) ,
\end{align}
which satisfies $ \Theta _{ l \infty } (s) = \Theta _{ l \infty } ( s + l ) $. 
This product converges absolutely for all $ s \in \mathbb{C} $. 
Hence we see that $ \Theta _{ l \infty } (s) $ is an entire function with zeros at 
$ s = 1/2 + l k + i r_{n} $ $ ( k \in \mathbb{Z} , n \geq 0 ) $. 
We note that the function $ \Theta _{ l \infty } (s) $ is also constructed 
by employing the idea of the zeta regularized product. (Cf. \cite{Kimoto}.) 

We are now in a position to state the functional equation of the higher Selberg zeta function. 

\begin{thm}[Main theorem]
Define the complete higher Selberg zeta function by
\begin{align*}
\Xi _{ l \infty } (s)
: = \Xi _{ l \infty , \mathrm{I} } (s) \cdot \Xi _{ l \infty , \mathrm{hyp} } (s) 
\cdot \Xi _{ l \infty , \mathrm{ell} } (s) \cdot \Xi _{ l \infty , \mathrm{par} } (s) . 
\end{align*}
Then $ \Xi _{ l \infty } (s) ^{-1} $ is an entire function of order $3$ with zeros at
$ s = 1/2 -l m \pm i r_{n} , ( m \geq 1 ,  n \geq 0 ) $ only, and satisfies
\begin{align}
& \Xi _{ l \infty } (s) = \frac { \Xi _{ l \infty } (s+l) }{ \Xi (s+l) } , \label{XSSL} \\
& \frac{ d^{2} }{ ds^{2} } \log \Xi _{ l \infty } (s) + \frac{ d^{2} }{ ds^{2} } \log \Xi _{ l \infty } (1-l-s) \label{DDXS} \\
= & \big( \frac{\pi}{l} \big)^{2} \sum_{ n=0 }^{ \infty } \Big\{ 
\sin ^{-2} \frac{ \pi  \big( s - \frac{1}{2}  - i r_{n} \big)}{l} 
+ \sin ^{-2} \frac{ \pi  \big( s - \frac{1}{2}  + i r_{n} \big) }{l} \Big\} . \notag
\end{align}
Also put $ \hat{ \Xi } _{ l \infty } (s) 
:= \Theta _{ l \infty } (s) \cdot \Xi _{ l \infty } (s) $, then $ \hat{ \Xi } _{ l \infty } (s) $ 
satisfies the functional equation:
\begin{align}
\hat{ \Xi } _{ l \infty } (s) \cdot \hat{ \Xi } _{ l \infty } (1-l-s) = C_{ l \infty } ,
\end{align}
where $ C_{ l \infty } $ is a non-zero constant.
\end{thm}

\begin{proof}
We first note that the relation (\ref{XSSL}) follows 
immediately from (\ref{HSZO}), (\ref{HIO}), (\ref{HEO}), (\ref{HPO}). 
This relation clearly shows that 
$ \Xi _{ l \infty } (s) ^{-1} $ is an entire function having zeros at $ s = 1/2 -l m \pm i r_{n} $. 
Next, by using (\ref{HC}), (\ref{IC}), (\ref{EC}) and (\ref{PC}), we have
\begin{align*}
\frac{ d^{2} }{ ds^{2} } \log \Xi _{ l \infty } (s) - \frac{1}{l} \frac{d}{ds} \log \Xi ( s + \frac{l}{2} ) 
= \sum_{ n=0 }^{ \infty } \big\{ \phi_{s} ( r_{n} ) + \phi_{s} ( - r_{n} ) \big\} .
\end{align*}
Replacing $ s $ by $ -s $ in the equation above, we have
\begin{align*}
\frac{ d^{2} }{ ds^{2} } \log \Xi _{ l \infty } (1-l-s) + \frac{1}{l} \frac{d}{ds} \log \Xi ( 1 - \frac{l}{2} -s ) 
= \sum_{ n=0 }^{ \infty } \big\{ \phi_{1-l-s} ( r_{n} ) + \phi_{1-l-s} ( - r_{n} ) \big\} .
\end{align*}
Here, from the equality
\begin{align*}
\sum_{ m = -\infty }^{ \infty } \frac{1}{ (z+lm)^{2} } = \big( \frac{ \pi }{l} \big) ^{2} \sin ^{-2} \frac{ \pi z }{l} ,
\end{align*}
we see that
\begin{align}
& \phi_{s} (r) + \phi_{s} (-r) + \phi_{ 1-l-s } (r) + \phi_{ 1-l-s } (-r) \label{phisr} \\
=& \sum_{ m= - \infty }^{ \infty } \frac{1}{ ( s+lm- \frac{1}{2} -ir )^{2} } 
+ \sum_{ m= - \infty }^{ \infty } \frac{1}{ ( s+lm- \frac{1}{2} +ir )^{2} } \notag \\
= & \big( \frac{\pi}{l} \big)^{2} \Big\{ 
\sin ^{-2} \frac{ \pi  \big( s - \frac{1}{2}  - i r \big)}{l} 
+ \sin ^{-2} \frac{ \pi  \big( s - \frac{1}{2}  + i r \big) }{l} \Big\} . \notag
\end{align}
Hence, combining these and the functional equation $ \Xi (s) = \Xi (1-s) $, we see that the relation (\ref{DDXS}) follows. 
Lastly, since  
\begin{align*}
\frac{ d^{2} }{ ds^{2} } \log \Theta _{ l \infty } (s) 
= - \big( \frac{\pi}{l} \big)^{2} \sum_{ n=0 }^{ \infty } \sin ^{-2} \frac{ \pi  \big( s - \frac{1}{2}  - i r_{n} \big)}{l} ,
\end{align*} 
it follows that
\begin{align*}
\frac{ d^{2} }{ ds^{2} } \log \hat{ \Xi } _{ l \infty }(s) \cdot \hat{ \Xi } _{ l \infty } (1-l-s) = 0 .
\end{align*} 
Since $ \hat{ \Xi } _{ l \infty }(s) \cdot \hat{ \Xi } _{ l \infty } (1-l-s) $ 
is invariant under the reflection $ s \leftrightarrow 1-l-s $, 
we conclude that it is a constant. 
This completes the proof.  
\end{proof}

By the expression of three factors 
$ \Xi _{ l \infty , \mathrm{I} } (s) $, $ \Xi _{ l \infty , \mathrm{ell} } (s) $, $ \Xi _{ l \infty , \mathrm{par} } (s) $, 
we can describe in principle all zeros and poles of the higher Selberg zeta function
$ \Xi _{ l \infty , \mathrm{hyp} } (s)
= \Xi _{ l \infty  } (s) \cdot \Xi _{ l \infty , \mathrm{I} } (s) ^{-1} \cdot 
\Xi _{ l \infty , \mathrm{ell} } (s) ^{-1} \cdot \Xi _{ l \infty , \mathrm{par} } (s) ^{-1} $. 
For example, when $ \varGamma = SL( 2 , \mathbb{Z} ) $ and $l=1$, 
zeros and poles of the higher Selberg zeta function 
$ \Xi _{ 1 \infty , \mathrm{hyp} } (s) $ are explicitly given as follows. \\

\noindent{\bf Zeros of $ \Xi _{ 1 \infty , \mathrm{hyp} } (s) $.}

(1) $ s= - 1/2 - k $ $( k \geq 0)$ ; \; order $ k+1 $. \\

\noindent{\bf Poles of $ \Xi _{ 1 \infty , \mathrm{hyp} } (s) $.}

(1) $ s=0 $ ; \quad order $1$,

(2) $ s = - 6k-j $ $ (k \geq 0 \; , j=2,3,4,5 ) $ ; \; order $ 6 k^{2} + 2jk + j - 1 $,

$ \;\; \quad  s= - 6k- j $ $ (k \geq 0 \; , j= 6,7 ) $ ; \; order $ 6 k^{2} + 2jk + j + 1 $,

(3) $ s = \rho / 2 - m $ ($m \geq 1$, \; $ \rho $: non-trivial zeros of $ \zeta (s) $) ,

(4) $ s = 1/2 -m \pm i r_{n} $ $ (m \geq 1 \; , n \geq 1 ) $.

\section{Higher Dirichlet $L$-functions}

We note that 
the higher Dirichlet $L$-function appears in the gamma factor of the higher Selberg zeta function for congruence subgroups.
In this section, we discuss the completion of the higher Dirichlet $L$-function itself by imitating the method of the previous section.
This is a generalization of the higher Riemann zeta function introduced in \cite{Cohen} and \cite{Matsuda}. 
In order to describe the functional equation, we use the explicit formula in place of the trace formula.

\subsection{Weil's Explicit Formula}

We first recall the functional equation of the Dirichlet $L$-function briefly. 
Let $ \chi $ be a primitive Dirichlet character modulo $ q \geq 2 $ or the trivial character $ \chi_{0} $ modulo $1$. 
We note that $ L ( s, \chi_{0} ) = \zeta (s) $ is the Riemann zeta function. 
Define the complete Dirichlet $L$-function $ \xi ( s, \chi ) $ by
\begin{align*}
\xi ( s , \chi ) 
:= \big\{ s(s-1) \big\} ^{ \delta_{ \chi } } \big( \frac{ \pi }{ q } \big) ^{ - s/2 } 
\Gamma \big( \frac{s+v}{2} \big) L ( s , \chi ) ,
\end{align*}
where $ \chi (-1) = (-1)^{v} $ $ (v = 0 \text{ or } 1) $, and $ \delta_{ \chi } $ is $1$ if $ \chi = \chi_{0} $ and $0$ otherwise. 
Then the function $ \xi ( s , \chi ) $ is an entire function of order $1$ and satisfies the functional equation 
$ \xi ( s , \chi ) = \epsilon ( \chi ) \xi (1-s, \bar{\chi} ) $, 
where $ \epsilon ( \chi ) $ is a constant modulus $1$ depending on $ \chi $ only. 
The Dirichlet $L$-function can be expressed both as an Euler product of prime numbers and as a Hadamard product of non-trivial zeros.
This fact leads the following formula due to Weil.

\begin{lem}[\cite{Weil}] 
Suppose that the function $ F : \mathbb{R} \to \mathbb{C} $ satisfies the following three conditions:

$(1)$ There exists $ a^{ \prime } > 0 $ such that the variation for $ F(x) e^{ ( \frac{1}{2} + a^{ \prime } ) |x| } $ on 
$ \mathbb{R}$ is bounded. 

$(2)$ $ F $ is normalized, that is
\begin{align*}
F(x) = \lim_{ h \to +0 } \frac{ F(x+h) + F(x-h) }{2},  \qquad \text{for all } x \in \mathbb{R} .
\end{align*}

$(3)$ There exists $ \epsilon > 0 $ such that
\begin{align*}
\frac{ F(x) + F(-x) }{2} = F(0) + O ( |x|^{ \epsilon } ), \qquad \text{for } x \to 0 .
\end{align*}

\noindent  For $ u > 0 $, put $ f(u) = F( - \log u ) $. Then we have
\begin{align*}
\lim_{ T \to \infty } \sum_{ | \mathrm{Im} \rho | < T } M ( f, \rho - \frac{1}{2} ) 
= & \delta _{ \chi } \big\{ M ( f, - \frac{1}{2} ) + M ( f, \frac{1}{2} ) \big\} - F(0) \log \frac{ \pi }{q} \\
+ & \frac{1}{ 2 \pi } \int_{ - \infty }^{ \infty } M (f, it) \;
\mathrm{Re} \Big\{ \frac{ \Gamma ^{ \prime } }{ \Gamma } \big( \frac{ 1+2v }{4} + \frac{it}{2} \big) \Big\} \; dt \notag \\
- & \sum_{ n=1 }^{ \infty } \sum_{ p : \mathrm{prime} } \frac{ \log p }{ p^{ \frac{n}{2} } }
\Big\{ \chi (p)^{n} F ( - \log p^{n} ) + \chi (p)^{-n} F ( \log p^{n} )  \Big\} . \notag
\end{align*}
Here $ M( f, s ) $ denotes the Mellin transform of $ f(u) $ as
\begin{align*}
M( f,s ) = \int_{0}^{ \infty } f(u) u^{s} \frac{du}{u}, 
\end{align*}
and $ \rho $ runs through the non-trivial zeros of the Dirichlet $L$-function $ L( s , \chi ) $ in $ 0 < \mathrm{Re} (s) < 1 $. \qed
\end{lem}

\subsection{Test Function and Gamma Factor}

By making use of the explicit formula, we investigate the higher Dirichlet $L$-functions for $ l \geq 1$: 
\begin{align*}
L_{ l \infty } ( s, \chi ) := \prod_{ m=1 }^{ \infty } L ( s+lm , \chi ) 
= \prod _{ m=1 } ^{ \infty } \prod _{p : \text{prime} } \big( 1 - \chi (p) \; p ^ { -s - lm } \big) ^{-1} , \quad
( \mathrm{Re} (s) > 1-l ) .
\end{align*}
It is seen that $ L_{ l \infty } ( s, \chi ) $ satisfies 
$ L_{ l \infty } ( s, \chi ) = L ( s+l , \chi ) \cdot L_{ l \infty } ( s+l , \chi ) $,
which gives an analytic continuation to the whole complex plane. 
We note that the logarithmic derivatives of $ L ( s , \chi ) $ and $ L_{ l \infty } ( s, \chi ) $ are given by 
\begin{align*}
\frac{d}{ds} \log L ( s , \chi ) 
= - \sum_{ n=1 }^{ \infty } \frac{ \chi (n) \Lambda (n) }{ n^{s} }
= - \sum_{ n=1 }^{ \infty } \sum_{ p: \mathrm{prime} } \frac{ \chi(p)^{n} \log p }{ p^{ \frac{n}{2} } }  p^{ - ( s - \frac{1}{2} ) n } ,
\end{align*}
and
\begin{align*}
\frac{ d^{2} }{ ds^{2} } \log L_{ l \infty } ( s, \chi )
= & \sum_{ n=1 }^{ \infty } \sum_{ p: \mathrm{prime} } \frac{ \chi (p)^{n} \log p }{ p^{ \frac{n}{2} } } 
\cdot \frac{ n \log p }{ p^{ \frac{ln}{2} } - p^{ - \frac{ln}{2} } } p^{ - ( s + \frac{l}{2} - \frac{1}{2} ) n } . 
\end{align*}
In view of these expressions, we choose the following test function: 
\begin{align*}
F(x) =  
\begin{cases} 
\displaystyle 
\frac{ x e^{ ( s + \frac{l}{2} - \frac{1}{2} )x }  }{ 2 \sinh \frac{lx}{2} } 
- \frac{1}{ l } e^{ ( s + \frac{l}{2} -\frac{1}{2}) x } 
& x < 0 \\
0 
& x \geq 0 , 
\end{cases} 
\end{align*}
for $ \mathrm{Re} (s) > 1-l $. 
This is essentially same the test function of (\ref{TestFun}) 
and satisfies the required condition of the explicit formula. 
Then the Mellin transform of $ f(u) $ is given by
\begin{align*}
M( f , it ) 
= \int_{0}^{ \infty } F( - \log u ) u^{ it } \frac { du }{u} 
= \int_{ - \infty }^{ \infty } F(-x) e^{itx} dx 
= \phi_{s} (-t) ,
\end{align*}
where 
\begin{align*}
\phi_{s} (t) := \sum_{ m=0 }^{ \infty } \Big\{ \frac{1}{ (s+lm+l-\frac{1}{2} +it )^{2} } 
- \frac{1}{l} \frac{1}{ s+lm + \frac{l}{2} - \frac{1}{2} +it }
+ \frac{1}{l} \frac{1}{ s+lm + \frac{ 3l }{2} - \frac{1}{2} +it } \Big\}. 
\end{align*}
Furthermore, we see that 
\begin{align*}
& - \sum_{ n=1 }^{ \infty } \sum_{ p : \mathrm{prime} } \frac{ \log p }{ p^{ \frac{n}{2} } }
\Big\{ \chi (p)^{n} F ( - \log p^{n} ) +  \chi (p)^{-n} F ( \log p^{n} )  \Big\} \\
= & - \frac{ d^{2} }{ ds^{2} } \log L_{ l \infty } ( s , \chi ) - \frac{1}{l} \frac{d}{ds} \log L \big( s+ \frac{l}{2} , \chi \big). 
\end{align*} 
Applying this test function, we calculate each terms of the explicit formula as follows. 

\begin{lem}
Retain the notation above. Then we have
\begin{align*}
& M ( f, - \frac{1}{2} ) + M ( f, \frac{1}{2} ) \\
= & - \frac{ d^{2} }{ ds^{2} } \log \Big\{ 
l^{ - \frac{2s}{l} } \cdot \Gamma \big( \frac{s+l}{l} \big) ^{-1} \Gamma \big( \frac{s+l-1}{l} \big) ^{-1} \Big\}
- \frac{1}{l} \frac{d}{ds} \log ( s + \frac{l}{2} ) ( s + \frac{l}{2} -1 ) , \notag
\end{align*}
and
\begin{align*}
- F(0) \log \frac{\pi}{q} = 0 
= - \frac{ d^{2} }{ ds^{2} } \log \big( \frac{\pi}{q} \big) ^{ \frac{s(s+l)}{4l} } 
- \frac{1}{l} \frac{d}{ds} \log \big( \frac{ \pi }{q} \big) ^{ - \frac{ s+\frac{l}{2} }{2} } .
\end{align*}
\end{lem}

\begin{proof}
We see that $ F (0) = 0 $ and 
\begin{align*}
& M ( f, - \frac{1}{2} ) + M ( f, \frac{1}{2} ) = \phi_{s} ( - \frac{i}{2} ) + \phi_{s} ( \frac{i}{2} ) \\
= & \sum_{ m=0 }^{ \infty } \Big\{ \frac{1}{ (s+lm+l)^{2} } \Big\}
- \frac{1}{l} \frac{1}{ s+ \frac{l}{2} } 
+ \sum_{ m=0 }^{ \infty } \Big\{ \frac{1}{ (s+lm+l-1)^{2} } \Big\}
- \frac{1}{l} \frac{1}{ s+ \frac{l}{2} - 1 } \\
= & \frac{ d^{2} }{ ds^{2} } \log \Gamma \big( \frac{s+l}{l} \big) \Gamma \big( \frac{s+l-1}{l} \big)
- \frac{1}{l} \frac{1}{ s+ \frac{l}{2} } - \frac{1}{l} \frac{1}{ s+ \frac{l}{2} -1 } .
\end{align*}
Hence the lemma holds. 
\end{proof}

\begin{lem}
We have
\begin{align*}
\frac{1}{ 2 \pi } \int_{ - \infty }^{ \infty } M (f, it) \;
\mathrm{Re} \Big\{ \frac{ \Gamma ^{ \prime } }{ \Gamma } \big( \frac{1+2v}{4} + \frac{it}{2} \big) \Big\} \; dt 
= - \frac{ d^{2} }{ ds^{2} } \log J_{l} (s+v) 
- \frac{1}{l} \frac{d}{ds} \log \Gamma \big( \frac {s+v + \frac{l}{2} }{2} \big) . 
\end{align*}
Here $ J_{l} (s) $ is defined by
\begin{align*}
J_{l} (s) 
:= k_{1} (l) ^{ - \frac{s}{l} } \cdot l^{ - \frac{ s(s+l) }{4l} } 
\cdot \prod_{ r=0 }^{ l-1 } \Gamma_{2} \big( \frac{s+2r+l}{2l} \big) \Gamma_{2} \big( \frac{s+2r+2l}{2l} \big) ,
\end{align*}
which satisfies 
$ J_{l} (s) / J_{l} (s-l) = \Gamma ( s/2 ) ^{-1} $.
\end{lem}

\begin{proof}
Using the residue theorem, we see that
\begin{align*}
& \frac{1}{ 2 \pi } \int_{ - \infty }^{ \infty } M (f, it) \;
\mathrm{Re} \Big\{ \frac{ \Gamma ^{ \prime } }{ \Gamma } \big( \frac{1+2v}{4} + \frac{it}{2} \big) \Big\} \; dt \\
= & \frac{1}{ 4 \pi } \int_{ - \infty }^{ \infty } 
\phi_{s}(t) \Big\{ \frac{ \Gamma ^{ \prime } }{ \Gamma } \big( \frac{1+2v}{4} + \frac{it}{2} \big) 
+ \frac{ \Gamma ^{ \prime } }{ \Gamma } \big( \frac{1+2v}{4} - \frac{it}{2} \big)\Big\} \; dt 
= - \sum_{ n=0 }^{ \infty } \phi_{s} \big( -\frac{4n+2v+1}{2} i \big) \\
= & - \sum_{ n=0 }^{ \infty }  \sum_{ m=0 }^{ \infty } \Big\{ \frac{1}{ ( s+v+lm+l+2n )^{2} } 
- \frac{1}{l} \frac{1}{ s+v+lm + \frac{l}{2} +2n }
+ \frac{1}{l} \frac{1}{ s+v+lm + \frac{ 3l }{2} + 2n } \Big\} .
\end{align*}
Putting $ n = l j_{1} + r_{1} $, $ m = 2 j_{2} + r_{2} $ and $ j_{1} + j_{2} = k $, we have
\begin{align*}
& - \sum_{ n=0 }^{ \infty }  \sum_{ m=0 }^{ \infty } 
\Big\{ \frac{1}{ ( s+v+lm+l+2n )^{2} } 
- \frac{1}{l} \frac{1}{ s+v+lm + \frac{l}{2} +2n }
+ \frac{1}{l} \frac{1}{ s+v+lm + \frac{ 3l }{2} + 2n } \Big\} \\
= & - \sum_{ r_{1} = 0 }^{ l-1 } \sum_{ r_{2} = 0 }^{1} \sum_{ j_{1} = 0 }^{ \infty } \sum_{ j_{2} = 0 }^{ \infty }
\Big\{ \frac{1}{ ( s+v+ 2l ( j_{1} + j_{2} ) +2 r_{1} + l r_{2} + l )^{2} } \\
& \qquad \qquad - \frac{1}{l} \frac{1}{  s+v+ 2l ( j_{1} + j_{2} ) +2 r_{1} + l r_{2} + \frac{l}{2} }
+ \frac{1}{l} \frac{1}{ s+v+ 2l ( j_{1} + j_{2} ) +2 r_{1} + l r_{2} + \frac{ 3l }{2}  } \Big\} \\
= & - \sum_{ r_{1} = 0 }^{ l-1 } \sum_{ k = 0 }^{ \infty } 
\Big\{ \frac{k+1}{ ( s+v+ 2lk +2 r_{1} + l )^{2} } + \frac{k+1}{ ( s+v+ 2lk +2 r_{1} + 2l )^{2} } \\
& \qquad \qquad - \frac{1}{l} \frac{k+1}{  s+v+ 2lk +2 r_{1} + \frac{l}{2} }
+ \frac{1}{l} \frac{k+1}{ s+v+ 2lk +2 r_{1} + \frac{ 5l }{2}  } \Big\} . 
\end{align*}
Since we have 
\begin{align*}
& \frac{ d^{2} }{ ds^{2} } \log \Gamma_{2} \big( \frac{ s+v+2r+l }{2l} \big)
= \frac{1}{ 4 l^{2} } ( 1 + \gamma ) 
+ \sum_{ k=0 }^{ \infty } \Big\{ \frac{ k+1 }{ (s+v+2lk+2r+l) ^{2} } - \frac{1}{ 4 l^{2} } \frac{1}{k+1} \Big\}, \\
& \frac{d}{ ds } \log \Gamma \big( \frac{ s+v+ 2r + \frac{l}{2} }{ 2l } \big)  
= - \frac{ \gamma }{2l}
- \sum_{ k=0}^{ \infty } \Big\{ \frac{1}{ s+v+ 2lk + 2r + \frac{l}{2} } - \frac{1}{2l} \frac{1}{k+1} \Big\} ,
\end{align*}
it follows that 
\begin{align*}
& \sum_{ k = 0 }^{ \infty } 
\Big\{ \frac{k+1}{ ( s+v+2lk+2r+l )^{2} } + \frac{k+1}{ ( s+v+2lk+2r+2l )^{2} } \\
& \quad \quad - \frac{1}{l} \frac{k+1}{  s+v+2lk+2r+ \frac{l}{2} }
+ \frac{1}{l} \frac{k+1}{ s+v+2lk+2r+ \frac{ 5l }{2}  } \Big\} \\
& = \frac{ d^{2} }{ ds^{2} } \log \Gamma_{2} \big( \frac{ s+v+2r+l }{ 2l } \big) \Gamma_{2} \big( \frac{ s+v+2r+2l }{ 2l } \big) 
- \frac{1}{l} \frac{d}{ ds } \log \Gamma \big( \frac{ s+v+2r + \frac{l}{2} }{ 2l } \big) .
\end{align*}
From Gauss-Legendre's multiplication formula (\ref{GaMult}), we conclude that 
\begin{align*}
& - \sum_{ r_{1} = 0 }^{ l-1 } \sum_{ k = 0 }^{ \infty } 
\Big\{ \frac{k+1}{ ( s+v+ 2lk +2 r_{1} + l )^{2} } + \frac{k+1}{ ( s+v+ 2lk +2 r_{1} + 2l )^{2} } \\
& \qquad \qquad \quad - \frac{1}{l} \frac{k+1}{  s+v+ 2lk +2 r_{1} + \frac{l}{2} }
+ \frac{1}{l} \frac{k+1}{ s+v+ 2lk +2 r_{1} + \frac{ 5l }{2}  } \Big\} \\
= &  - \sum_{ r = 0 }^{ l-1 } \Big\{ 
\frac{ d^{2} }{ ds^{2} } \log \Gamma_{2} \big( \frac{ s+v+2r+l }{ 2l } \big) \Gamma_{2} \big( \frac{ s+v+2r+2l }{ 2l } \big) 
+ \frac{1}{l} \frac{d}{ ds } \log \Gamma \big( \frac{ s+v+2r + \frac{l}{2} }{ 2l } \big) \Big\} \\ 
= & - \frac{ d^{2} }{ ds^{2} } \log J_{l} (s+v) 
- \frac{1}{l} \frac{d}{ds} \log \Gamma \big( \frac {s+v + \frac{l}{2} }{2} \big) .
\end{align*}
This complete the proof. 
\end{proof}

\subsection{Complete Higher Dirichlet $L$-Function}

To describe the functional equation of $ L _{ l \infty } ( s , \chi ) $, 
we introduce the function
\begin{align*}
\theta _{ l \infty } ( s , \chi )
:= \prod_{ n=1 }^{ \infty }
\big( 1 - e^{ \frac{ 2 \pi i }{l} ( \rho _{n} - s ) } \big) .
\end{align*}
Here $ \rho _{n} = 1/2 + i t_{n} $ $ ( \mathrm{Re} ( t_{n} ) \geq 0 ) $ 
is the non-trivial zero of the Dirichlet $L$-function $ L ( s, \chi ) $ on the upper half plane. 
It is seen that $ \theta _{ l \infty } ( s, \chi ) $ is an entire function of order $2$ with zeros at 
$ s = 1/2 + i t_{n} + l k $, $ ( n \geq 1 , \; k \in \mathbb{Z} ) $, 
and satisfies 
$ \theta _{ l \infty } ( s , \chi ) = \theta _{ l \infty } ( s + l , \chi ) $. 
Now, combining above lemmas and applying the method of Theorem 3.6.1, 
we obtain the functional equation of the higher Dirichlet $L$-function as follows.
\begin{thm}
Define the complete higher Dirichlet $L$-function by
\begin{align*}
\xi _{ l \infty } ( s , \chi )
: = & \Big\{ \exp \big( - \frac{ 2s }{l} \log l \big) 
\cdot \Gamma \big( \frac{s+l}{l} \big)^{-1} \Gamma \big( \frac{s+l-1}{l} \big)^{-1} \Big\} ^{ \delta_{ \chi } } \\
\times & \exp \big\{ \frac{ s(s+l) }{ 4l } \log \big( \frac{ \pi }{q} \big) 
- \frac{s+v}{l} \log \big( k_{1} (l) \big) - \frac{ (s+v)(s+v+l) }{ 4l } \log l \big\} \\
\times & \prod_{ r=0 }^{ l-1 } \Gamma_{2} \big( \frac{s+v+2r+l}{2l} \big) \Gamma_{2} \big( \frac{s+v+2r+2l}{2l} \big) 
\times L_{ l \infty } ( s, \chi ).
\end{align*}
Then $ \xi _{ l \infty } ( s , \chi ) $ is an entire function of order $2$ with zeros at $ s = \rho -lm $ $ ( m \geq 1 )$, 
and satisfies 
$ \xi _{ l \infty } ( s , \chi ) = \xi ( s+l, \chi ) \cdot \xi _{ l \infty } (s+l ,\chi ) $. 
In addition, 
$ \hat{ \xi } _{ l \infty } ( s , \chi ) :=  \theta _{ l \infty } (s , \chi ) ^{-1} \cdot \xi _{ l \infty } ( s, \chi ) $ 
satisfies the functional equation:
\begin{align*}
\hat{ \xi } _{ l \infty } (s , \chi ) \cdot \hat{ \xi } _{ l \infty } (1-l-s, \bar{ \chi } ) 
= \exp \big\{ a_{ l \infty } ( \chi ) s + b_{ l \infty } ( \chi ) \big\} ,
\end{align*}
where $ a_{ l \infty } ( \chi ) $ and $ b_{ l \infty } ( \chi ) $ are constants. \qed
\end{thm}


\end{document}